\documentclass[12pt]{amsart}
\usepackage{amsfonts, amssymb, amsmath, amsthm, latexsym, array}
\usepackage{fullpage,color}
\usepackage{verbatim,euscript}

\numberwithin{equation}{section}

\usepackage[colorlinks=true,linkcolor=blue,urlcolor=violet,citecolor=magenta]{hyperref}

\input {cyracc.def}
 \font\tencyr=wncyr10 
\font\tencyi=wncyi10 
\font\tencysc=wncysc10 
\def\rus{\tencyr\cyracc}
\def\rusi{\tencyi\cyracc}
\def\rusc{\tencysc\cyracc}

\newtheorem{thm}{Theorem}[section] 
\newtheorem{conj}[thm]{Conjecture}
\newtheorem{lm}[thm]{Lemma}
\newtheorem{cor}[thm]{Corollary}
\newtheorem{prop}[thm]{Proposition}

\theoremstyle{remark}
\newtheorem{rmk}[thm]{Remark}

\theoremstyle{definition}
\newtheorem{ex}[thm]{Example}
\newtheorem*{ex-bn}{Example}

\newcommand {\be}{{\mathfrak b}}
\newcommand {\ce}{{\mathfrak c}}

\newcommand {\g}{{\mathfrak g}}

\newcommand {\me}{{\mathfrak m}}
\newcommand {\n}{{\mathfrak n}}
\newcommand {\fN}{\mathfrak{N}}
\newcommand {\p}{{\mathfrak p}}

\newcommand {\te}{{\mathfrak t}}
\newcommand {\ut}{{\mathfrak u}}

\newcommand {\cf}{{\mathcal F}}
\newcommand {\cl}{{\mathcal L}}
\newcommand {\cM}{{\mathcal M}}

\newcommand {\co}{{\mathcal O}}
\newcommand {\cs}{{\mathcal S}}

\newcommand {\BZ}{{\mathbb Z}}
\newcommand {\BN}{{\mathbb N}}
\newcommand {\BQ}{{\mathbb Q}}

\newcommand {\BC}{{\mathbb C}}

\newcommand {\isom}{\stackrel{\sim}{\rightarrow}}

\newcommand {\esi}{\varepsilon}
\newcommand{\lb}{\lambda}
\newcommand{\ap}{\alpha}
\newcommand{\vp}{\varphi}

\renewcommand{\le}{\leqslant}
\renewcommand{\ge}{\geqslant}
\newcommand{\curge}{\succcurlyeq}
\newcommand{\curle}{\preccurlyeq}

\newcommand{\Xp}{\mathfrak X_+}

\newcommand{\eus}{\EuScript}

\newcommand {\hot}{{\mathsf{ht}}}

\newcommand {\Lie}{{\mathsf{Lie}}}

\newcommand {\Hom}{\operatorname{Hom}}

\newcommand {\rk}{{\mathsf{rk}}}

\newcommand {\GR}[2]{{\textrm{{\bf #1}}}_{#2}}

\newcommand {\ov}{\overline}
\newcommand {\un}{\underline}
\newcommand {\bD}{{\boldsymbol\Delta}}
\newcommand {\bU}{{\boldsymbol U}}
\newcommand {\bX}{{\boldsymbol X}}
\newcommand {\bZ}{{\boldsymbol Z}}
\newcommand {\vbt}{V_{\bar\theta}}

\newfam\eusfam%
\font\Bbbfont=msbm10 scaled 1200%
\def\bbk{\hbox {\Bbbfont\char'174}}

\begin{document}
\setlength{\parskip}{3pt plus 5pt minus 0pt}
\hfill { {\scriptsize April 23, 2009}}
\vskip1ex

\title[Generalised Kostka-Foulkes  polynomials]
{Generalised Kostka-Foulkes  polynomials and cohomology of line bundles on
homogeneous vector bundles}
\author[D.\,Panyushev]{Dmitri I.~Panyushev}
\address[]{Independent University of Moscow,
Bol'shoi Vlasevskii per. 11, 119002 Moscow, \ Russia
\hfil\break\indent
Institute for Information Transmission Problems, B. Karetnyi per. 19, Moscow 127994
}
\email{panyush@mccme.ru}
\maketitle

\tableofcontents
\section*{Introduction}
\noindent
Let $G$ be a  semisimple algebraic group with Lie algebra $\g$. 
We consider 
generalisations of Lusztig's  $q$-analogue of weight multiplicity.  
Fix a maximal torus $T\subset G$. Let $m_\lb^\mu$ be the multiplicity of weight 
$\mu$ in a simple $G$-module $V_\lb$ with highest weight $\lb$.  
Lusztig's $q$-analogues $\me_\lb^\mu(q)$ (also known as Kostka-Foulkes polynomials 
for the root system of $G$)
are certain polynomials in $q$ such that $\me_\lb^\mu(1)=m_\lb^\mu$.
 A recent survey of their properties, with an eye towards combinatorics, is given in \cite{n-r}.
These polynomials arise in numerous problems of representation theory, geometry, and
combinatorics. Work of Lusztig \cite{lus} and Kato \cite{kato} shows that, for
$\lb$ and $\mu$ dominant,  $\me_\lb^\mu(q)$ are connected with certain Kazhdan-Lusztig polynomials for the affine Weyl group associated with $G$.
To define $\me_\lb^\mu(q)$, one
first considers a $q$-analogue of Kostant's partition function, $\eus P$. It is conceivable to 
replace the set of positive roots, $\Delta^+$, occurring in the definition of $\eus P$ 
with an arbitrary 
finite multiset $\Psi$ in the character group $\mathfrak X$ of  $T$. 
If the elements of $\Psi$ belong
to an open half-space of $\mathfrak X\otimes \BQ$ (this is our {\sl first\/} hypothesis on $\Psi$), 
then we still obtain certain polynomials
$\me_{\lb,\Psi}^\mu(q)$.  
We always  assume that $\lb$ is dominant, whereas $\mu\in\mathfrak X$
can be arbitrary. In this article, we are interested in the non-negativity problem
for the coefficients of
$\me_{\lb,\Psi}^\mu(q)$.
For Lusztig's  $q$-analogues, this problem 
has been considered by Broer. He proved that  
%
 $\me_{\lb}^{\mu}(q)$  has non-negative coefficients for any  $\lb\in\Xp$
if and only if $(\mu,\ap^\vee)\ge -1$ for all $\ap\in\Delta^+$
(see \cite[Theorem\,2.4]{br93} and \cite[Prop.\,2(iii)]{br97}).

%

Our first goal is to provide sufficient conditions for $\me_{\lb,\Psi}^\mu(q)$ to have non-negative coefficients.
Let $B$ be the Borel subgroup of $G$ corresponding to $\Delta^+$
(i.e., the roots of $B$ are positive!) and $\Xp$ the set of dominant weights.
The {\sl second\/} hypothesis is that $\Psi$ is assumed to be the
multiset of weights for a $B$-submodule  $N$ of a $G$-module $V$.
Then  $\me_{\lb,\Psi}^\mu(q)$ is said to be a 
{\it generalised Kostka-Foulkes polynomial}.
Let $P\supset B$ be any parabolic subgroup normalising $N$ and $G\times_P N$ the 
corresponding homogeneous vector bundle on $G/P$. We obtain a relation 
between the Euler characteristic of induced line bundles $\cl$ on the
$G\times_P N$ and generalised Kostka-Foulkes polynomials. Using the collapsing
$G\times_P N\to G{\cdot}N\subset V$, we get a vanishing result for 
$H^i(G\times_P N,\cl)$, $i\ge 1$, and conclude
that  $\me_{\lb,\Psi}^\mu(q)$ has non-negative
coefficients for all $\lb\in\Xp$ if $\mu$ is sufficiently large.
An explicit lower bound for $\mu$ is also given, see Section~\ref{sect:cohomol}.
This approach is based on the Grauert-Riemenschneider vanishing theorem. 
We also notice that Broer's formula for $\frac{d}{dq}\me_\lb^\mu(q)$ \cite{br95} can be generalised to $\me_{\lb,\Psi}^\mu(q)$. 
The most natural examples of generalised Kostka-Foulkes polynomials occur if $\Psi\subset
\Delta^+$. For instance, one can take $N$ to be a $B$-stable ideal in $\Lie (B,B)\subset \g$.

Our second goal is to study in details the special case in which $\Psi=\Delta^+_s$,
the set of short positive roots. The required $B$-submodule, $\vbt^+$, lies in $\vbt$, where 
$\bar\theta$ is the short dominant root. 
The polynomials $\ov{\me}_\lb^\mu(q):=
\me_{\lb, \Delta^+_s}^\mu(q)$ are said to be {\it short $q$-analogues}. 
The numbers $\ov{\me}_\lb^\mu(1)$ appeared already in work of Heckman \cite{heck},
and a geometric interpretation of  $\ov{\me}_\lb^0(q)$ given in \cite{ww}
shows that $\ov{\me}_\lb^0(q)$ have non-negative coefficients. 
Let $\Delta^+_l$ be the set of long positive roots, $W_l$  the (normal) 
subgroup of $W$ generated by all $s_\ap$ ($\ap\in\Delta^+_l$),  and
$\rho_l$  the half-sum of the long positive roots.
Approach of Section~\ref{sect:cohomol} enables us to prove that 
$\ov{\me}_\lb^\mu(q)$ has nonnegative coefficients whenever  $\mu+\rho_l\in\Xp$
(Cor.~\ref{cor:rho-l}).
But to obtain  exhaustive results, we take  another path.
We consider the {\it shifted\/} (=\,dot) action of $W_l$ on $\mathfrak X$, 
$(w,\mu)\mapsto  w\odot\mu=w(\mu+\rho_l)-\rho_l$, and show that 
$\ov{\me}_\lb^{w\odot\mu}(q)=(-1)^{\ell(w)}\ov{\me}_\lb^\mu(q)$. 
Therefore $\ov{\me}_\lb^\mu(q)\equiv 0$ if $\mu$ is not regular 
relative to the shifted $W_l$-action, and it suffices to
consider $\ov{\me}_\lb^\mu(q)$ only for  $\mu$ that are dominant with respect to
$\Delta^+_l$. For a $\Delta^+_l$-dominant $\mu$,  we prove that 
$\ov{\me}_\lb^\mu(q)$ has non-negative coefficients for all $\lb\in\Xp$ 
if and only if  $(\mu,\ap^\vee)\ge -1$ for all $\ap\in\Delta^+_s$, see Theorem~\ref{thm:main}.
This is an extension of Broer's results in \cite[Sect.\,2]{br93}. Again, this stems from a
careful study of cohomology of line bundles on $G\times_B \vbt^+$.
In these considerations, it is important that  $W$ is a semi-direct product 
$W(\Pi_s)\ltimes W_l$, where the first group
is generated by the short simple reflections.
Modifying approach of R.~Gupta  \cite{rkg87}, 
we define analogues of Hall-Littlewood polynomials (Section~\ref{sect:HL}). These polynomials in $q$, denoted $\ov{P}_\lb(q)$,  are indexed by 
$\lb\in\Xp$ and form a $\BZ$-basis for the $q$-extended character ring $\boldsymbol{\Lambda}[q]$ of $G$. Let $\chi_\lb$ be the character of $V_\lb$ and $H$
the connected semisimple subgroup of $G$ whose root system is $\Delta_l$.
The polynomials $\ov{P}_\lb(q)$ interpolate between $\chi_\lb$
(at $q=0$) and a certain sum of irreducible characters of  $H$ (at $q=1$).
We obtain some orthogonality relations for $\ov{P}_\lb(q)$ and
show that  $\chi_\lb=\sum_{\mu\in\Xp}
\ov{\me}_{\lb}^{\mu}(q) \ov{P}_\mu(q)$.
Moreover,  
the whole theory developed by R.~Gupta  in  \cite{rkg87, rkg87-A}
can be extended to this setting.
For instance, we prove a version of Kato's  identity \cite[1.3]{kato} and point out a scalar
product in  $\boldsymbol{\Lambda}[q]$ such that $\{\ov{P}_\lb(q)\}_{\lb\in\Xp}$ to be an orthogonal basis. In a sense, the reason for such an extension is that $G{\cdot}\vbt=:\fN(\vbt)$ is the null-cone in $\vbt$, and, as well as the nilpotent cone $\fN\subset\g$, 
this variety is an irreducible normal complete intersection. 
On the other hand, Theorem~\ref{thm:main} yields vanishing of higher cohomology of the
structure sheaf $\co_{G\times_B \vbt^+}$, and, together with
\cite{kov00}, this implies that $\fN(\vbt)$ has only rational singularities.

We conjecture that if  $\mu$ satisfies vanishing conditions of Theorem~\ref{thm:main},  
then $\ov{\me}_\lb^\mu(q)$ can be interpreted as the "jump polynomial"  associated with a filtration of a subspace of $V_\lb^\mu$, see Subsection~\ref{subs:jump}. 
This is inspired by \cite{rkb89}.

\vskip1ex
{\small {\bf Acknowledgements.}  
This work was completed during my stay at I.H.\'E.S.
(Bures-sur-Yvette) in Spring 2009. 
I am grateful to this institution for the warm hospitality and support.}

\section{Notation} 
\label{sect:notation}

Let $G$ be a connected semisimple algebraic group of rank $r$, with a fixed Borel subgroup
$B$ and a maximal torus $T\subset B$.
The corresponding triangular decomposition of $\g=\Lie(G)$ is
$\g=\ut^-\oplus\te\oplus\ut$ and $\be=\te\oplus\ut$.
The character group of $T$ is denoted by $\mathfrak X$. 
Let $\Delta$ be the root system of $(G,T)$. Then $B$ determines the set of positive roots
$\Delta^+$ and the monoid of dominant weights $\Xp$.

\textbullet \ \ $\Pi$ 
is the set of simple roots in $\Delta^+$;
\\ \indent
\textbullet \ \ $\vp_1,\dots,\vp_r$ are the fundamental weights in $\Xp$.

\noindent
Write $W$ for the Weyl group and $s_\ap$ for the reflection corresponding to $\ap\in\Delta^+$. Set $\mathsf{N}(w)=\{\ap\in\Delta^+\mid w\ap \in -\Delta^+\}$ and 
$\esi(w)=(-1)^{\ell(w)}$, where $\ell(w)=\# \mathsf{N}(w)$ is the usual
length function on $W$. For  $\mu\in\mathfrak X$, let $\mu^+$ denote the unique dominant
element in $W\mu$.
We fix a $W$-invariant scalar product $(\ ,\ )$ on 
$\mathfrak X\otimes_\BZ \BQ$. As usual, $\ap^\vee=2\ap/(\ap,\ap)$ for
$\ap\in\Delta$.  
For any $\lb\in\Xp$, we choose a simple highest weight module $V_\lb$;
$V_\lb^\mu$ is the $\mu$-weight space in $V_\lb$ and
$m_{\lb}^{\mu}=\dim V_\lb^\mu$.

We consider two partial orders in $\mathfrak X$. For $\mu,\nu\in\mathfrak X$,
\begin{itemize}
\item the {\it root order}  is defined by letting
$\mu\curle \nu$ if and only if  $\nu-\mu$ lies in the monoid generated by $\Delta^+$;
notation $\mu\prec \nu$ means that $\mu\curle \nu$ and $\mu\ne\nu$;
\item the {\it dominant order} is defined by letting
$\mu\lessdot \nu$ if and only if  $\nu-\mu\in\Xp$. 
\end{itemize}

If $\Psi$ is a finite multiset in $\mathfrak X$, then $|\Psi|$ is the sum of all elements of
$\Psi$ (with respective multiplicities). Recall that $|\Delta^+|/2 =\vp_1+\ldots+\vp_r$, and this 
quantitiy is denoted by $\rho$.

Let $P$ be a parabolic subgroup of $G$. For a $P$-module $N$, let 
$G\times_P N$ denote the homogeneous $G$-vector bundle on $G/P$ whose fibre 
over $\{P\}\in G/P$ is $N$;
we write $\cl_{G/P}(V)$ for the locally free $\co_{G/P}$-module of its sections. If $N$ is a submodule of a $G$-module, then
the natural morphism $f: G\times_P N\to G{\cdot}N$ is projective and $G$-equivariant.
It is a {\it collapsing\/} in the sense of Kempf \cite{kempf}. 
Recall that $G{\cdot}N$ is a closed subvariety of $V$, since $N$ is $P$-stable.
If $\dim G\times_P N=\dim G{\cdot}N$,
then $f$ is said to be {\it generically finite}.
If $N'$ is another $P$-module, then $G\times_P (N\oplus N')$ is a vector bundle
on $G\times_P N$ with sheaf of sections $\cl_{G\times_P N}(N')$.

For any graded $G$-module $\eus C=\oplus_j \eus C_j$ with $\dim \eus C_j< \infty$,
its $G$-{\it Hilbert series\/} is defined by
\[
   \mathcal H_G(\eus C; q)=\sum_j \sum_{\lb\in\Xp} \dim \Hom_G(V_\lb, \eus C_j)
     e^\lb q^j \in \BZ [\mathfrak X] [[q]] .
\]

\section{Main definitions and first properties} 
\label{sect:result}

\noindent
Let $V$ be a finite-dimensional rational $G$-module and  $N$ a $P$-stable subspace 
of $V$.
We assume that the $T$-weights
occurring in $N$ lie in an open half-space of $\mathfrak X\otimes_\BZ \BQ$.
(This hypothesis implies that all $v\in N$ are unstable vectors in the sense of Geometric Invariant Theory.)
Counting each $T$-weight according to its multiplicity in $N$, we get a finite multiset
$\Psi$ in $\mathfrak X$.
The {\it generalised partition function}, $\eus P_{\Psi}$, is defined by the series
      $\displaystyle \frac{1}{\prod_{\ap\in\Psi} (1-e^\ap)}=\sum_{\nu} \eus P_{\Psi}(\nu)e^\nu$.
Accordingly, its $q$-analogue is defined by 
\[
      \frac{1}{\prod_{\ap\in\Psi} (1-qe^\ap)}=\sum_{\nu} \eus P_{\Psi,q}(\nu)e^\nu .
\] 
In view of our assumption on $N$, the numbers $\eus P_{\Psi}(\nu)$ are well-defined,
and $\eus P_{\Psi,q}(\nu)$ is a polynomial in $q$, with non-negative integer coefficients.
Clearly, $\eus P_{\Psi,q}(\nu)$ counts the "graded occurrences" of $\nu$ in the symmetric
algebra  $\mathcal S^\bullet(N)$. That is, $[q^j]\eus P_{\Psi,q}(\nu)=\dim \mathcal (S^jN)^\nu$.

For $\lb\in\Xp$ and $\mu\in \mathfrak X$, define the polynomials $\me^\mu_{\lb,\Psi}(q)$ by
\begin{equation}   \label{main-def}
   \me^\mu_{\lb,\Psi}(q)=\sum_{w\in W} \esi(w) \eus P_{\Psi,q}(w(\lb+\rho)-(\mu+\rho)) .
\end{equation}
This definition makes sense for any  multiset $\Psi$. But we require that
our $\Psi$ to be always the multiset of weights of a $P$-submodule of a $G$-module,
since we are going to exploit geometric methods.

For $N=\ut\subset\g$ and $\Psi=\Delta^+$, one obtains Lusztig's $q$-analogues of weight multiplicity \cite{lus}
(= Kostka-Foulkes polynomials for $\Delta$), and $\me_{\lb,\Delta^+}^\mu(1)=m_\lb^\mu$.
 Therefore, $\me^\mu_{\lb,\Psi}(q)$ is said to be
a  $(\Psi,q)$-{\it analogue of weight multiplicity\/} or {\it generalised 
Kostka-Foulkes polynomial}. 
If $\Psi=\Delta^+$, we will omit the subscript $\Delta^+$ in previous formulae. 

As $\me_{\lb,\Psi}^\mu(q)$ is a polynomial in $q$, one might be interested in its derivative.
For $\Psi=\Delta^+$, a nice formula for $\frac{d}{dq}\me_\lb^\mu(q)$ is found by Broer 
\cite[p.\,394]{br95}.
We notice that his method works in general, and it is more natural to begin with
a formula for the derivative of $\eus P_{\Psi,q}(\nu)$.

\begin{thm}
$\displaystyle \frac{d}{dq}\eus P_{\Psi,q}(\nu)=
\sum_{\gamma\in\Psi}\sum_{n\ge 1}q^{n-1}\eus P_{\Psi,q}(\nu-n\gamma)$.
\end{thm}
\begin{proof}
The derivative $\frac{d}{dq}\eus P_{\Psi,q}(\nu)$ equals the coefficient of $t$ in
the  expansion of  $\eus P_{\Psi,q+t}(\nu)$.  Let the polynomials $\eus R_{n,\mu}(q)$ be
defined by the generating function
\[
   \prod_{\ap\in\Psi} \frac{1-qe^\ap}{1-(q+t)e^\ap}= 
   \frac{\sum_{\nu} \eus P_{\Psi,q+t}(\nu)e^\nu}{\sum_{\nu} \eus P_{\Psi,q}(\nu)e^\nu}=:
   \sum_\mu \sum_{n\ge 0} \eus R_{n,\mu}(q)e^\mu t^n .
\] 
It is easy to compute these polynomials
for $n=0,1$.  First, taking $t=0$, we obtain  $\sum_\mu \eus R_{0,\mu}(q)e^\mu=1$.
Second, we have
\[
  \sum_\mu \eus R_{1,\mu}(q) e^\mu=\left[  \prod_{\ap\in\Psi} \frac{1-qe^\ap}{1-(q+t)e^\ap}\right]^\prime_t \vert_{t=0}=\sum_{\ap\in\Psi}\frac{e^\ap}{1-qe^\ap}=\sum_{\ap\in\Psi}
  \sum_{n\ge1} q^{n-1}e^{n\ap} .
\]
Hence $\eus R_{1,\mu}(q)=\begin{cases}  q^{n-1} & \text{if $\mu=n\ap$}, \ap\in\Psi \\
0, & \text{otherwise}.   \end{cases}$ 
\\  
Next, $\sum_\nu \eus P_{\Psi,q+t}(\nu)e^\nu=\sum_{n,\mu,\gamma} \eus R_{n,\mu}(q)
\eus P_{\Psi,q}(\gamma) e^{\mu+\gamma}t^n$. 
Hence 
\[
\eus P_{\Psi,q+t}(\nu)e^\nu= \sum_{n,\mu}  \eus R_{n,\mu}(q)\eus P_{\Psi,q}(\nu-\mu)t^n ,
\] 
 and extracting the coefficient of $t$ we get the assertion.
\end{proof}

\begin{cor}
$\displaystyle \frac{d}{dq}\me_{\lb,\Psi}^\mu(q)=\sum_{\gamma\in\Psi}\sum_{n\ge 1}q^{n-1}
\me_{\lb,\Psi}^{\mu+n\gamma}(q)$.
\end{cor}

It would be nice to have a formula for the degree of these polynomials
and necessary conditions for $\me_{\lb,\Psi}^\mu(q)$ to be nonzero. 
For Lusztig's $q$-analogues, it is easily seen that $\me_{\lb}^\mu(q)\ne 0$ if and only
if $\mu\curle \lb$, and $\deg \me_{\lb}^\mu(q)= \hot(\lb-\mu)$.
However, if  $\Psi$ is arbitrary, i.e., there is no relation between $\Delta^+$ and
$\Psi$, then  it is impossible to compare the degrees
of different summands in Equation~\eqref{main-def}. 
The only general assertion we can prove concerns the 
case in which $\Psi\subset \Delta^+$.

\begin{lm}
Suppose that\/ $\Psi\subset \Delta^+$. Then $\me_{\lb,\Psi}^\lb(q)=1$ 
and  if\/ $\me_{\lb,\Psi}^\mu(q) \ne 0$, then $\mu\curle \lb$.
\end{lm}
\noindent
Note that 
if $\me_{\lb,\Psi}^\mu(q) \ne 0$, then it is not necessarily true that $\lb-\mu$
lies in the monoid generated by $\Psi$.

\section{Cohomology of line bundles 
and  generalised Kostka-Foulkes polynomials}
\label{sect:cohomol}

\subsection{Statement of main results}
We assume that $P\supset B$ and choose a Levi subgroup $L\subset P$ such that 
$L\supset T$.  Write $\n$ for the nilpotent radical of $\p=\Lie(P)$, and $\Delta(\n)$ for the
roots of $\n$; hence $\Delta(\n)\subset \Delta^+$.
Let $\mathfrak X^P$ denote the character group of $P$. Obviously, $\mathfrak X^P$ is the character group of the central torus in $L$, and we may identify $\mathfrak X^P$ with
a subgroup of $\mathfrak X$.
Then $\Xp^P=\Xp\cap \mathfrak X^P$ is the monoid of $P$-dominant weights, 
i.e., the dominant weights $\lb$ such that $P$ stabilises a nonzero line in
$V_\lb$.  Let $\rho_P$ be the sum of those fundamental weights that belong to $\Xp^P$.

In this section, we prove the following two theorems:

\begin{thm}  \label{rez-2}   Set $\bZ=G\times_P N$. For $\mu\in \mathfrak X^P$, let 
$\cl_\bZ(\mu)^\star$ be the dual of the sheaf of sections of the line bundle
$G\times_P (N\oplus \BC_\mu)\to \bZ$. Then
\begin{itemize}
\item[\sf (i)]  \ 
$H^i(\bZ,\cl_\bZ(\mu)^\star)=0$ for all $i\ge 1$
whenever  $\mu \gtrdot \rho_P+|\Psi|-|\Delta(\n)|$.
\item[\sf (ii)]  \ If the collapsing 
$\bZ\to G{\cdot}N$ is generically finite,  
then $H^i(\bZ,\cl_\bZ(\mu)^\star)=0$ for all $i\ge 1$
whenever  $\mu \gtrdot |\Psi|-|\Delta(\n)|$.
\end{itemize}
\end{thm}

\begin{thm}  \label{rez-1}  
Suppose $N$ is $P$-stable and $\mu\in \mathfrak X^P$.
\begin{itemize} 
\item[\sf (i)]  \ If $\mu \gtrdot \rho_P+ |\Psi|- |\Delta(\n)|$, then 
$\me^\mu_{\lb,\Psi}(q)$ has non-negative coefficients for any $\lb\in\Xp$.
\item[\sf (ii)]  If the collapsing 
$G\times_P N\to G{\cdot}N$ is generically finite, then 
$\me^\mu_{\lb,\Psi}(q)$ has non-negative coefficients for any $\lb\in\Xp$
whenever $\mu \gtrdot  |\Psi|- |\Delta(\n)|$.
\end{itemize}
\end{thm}

\noindent
(Note that $|\Psi|,|\Delta(\n)|\in \mathfrak X^P$. Hence both inequalities concern weights
lying in $\mathfrak X^P$.)

\noindent
Actually, Theorem~\ref{rez-1} follows from Theorem~\ref{rez-2} and a relation between
$(\Psi,q)$-analogues and cohomology of line bundles, see Theorem~\ref{thm:svyaz-P}  below. 
Such an approach to $(\Psi,q)$-analogues is inspired by work of Broer \cite{br93,br94}.

\subsection{Algebraic-geometric facts}
For future reference, we recall some standard results in the form that we need below.
Let $U$ be the total space of a line bundle on an algebraic variety $Z$ and
$\pi: U\to Z$ be the corresponding projection. If $\mathcal E$ is a locally free $\co_Z$-module,
then $\mathcal E^\star$ is its dual.

\begin{lm}  \label{lm:odin}   
Let $\mathcal F$ be the sheaf  of sections of $\pi$.
\begin{itemize}
\item[\sf (i)] \  
If $\cl$ is a locally free $\co_Z$-module of finite type, then
$\pi_\ast(\pi^\ast\cl)=\bigoplus_{n\ge 0}(\cl\otimes(\mathcal F^{\otimes n})^\star)$.
\item[\sf (ii)] \ If $\mathcal G$ is a quasi-coherent sheaf on $U$, then
$H^i(U,\mathcal G)=H^i(Z, \pi_\ast\mathcal G)$ for all $i$.
\end{itemize}
\end{lm}
\begin{proof}
(i) Use the "projection formula" and the equality 
$\pi_\ast(\co_U)=\bigoplus_{n\ge 0}(\mathcal F^{\otimes n})^\star$.
\\
(ii)  This is true because $\pi$ is an affine morphism.
\end{proof}

\noindent
Thus, vanishing of higher cohomology for $\pi^\ast\cl$ will imply that
for $\cl\otimes(\mathcal F^{\otimes n})^\star$ for all ${n\ge 0}$.
The following is a special case of the Grauert--Riemenschneider theorem in Kempf's 
version (\cite[Theorem\,4]{kempf}):

\begin{thm}    \label{thm-GR}
Let $\omega_U$ denote the canonical bundle on $U$.
Suppose there is   a proper generically finite morphism $U\to X$ onto an affine variety $X$.
Then $H^i(U,\omega_U)=0$ for all $i\ge 1$.
\end{thm}

\subsection{Proof of Theorem~\ref{rez-2}}
Recall that $N$ is a $P$-submodule of a $G$-module $V$, $\Psi$ is the corresponding multiset of weights, and $\Psi$ belongs to an open half-space of 
$\mathfrak X\otimes_\BZ\BQ$.
Our goal is to obtain a sufficient condition for vanishing of higher cohomology of 
line bundles on $\bZ:=G\times_PN$.

For $\mu\in \mathfrak X^P_+$, let  $\BC_\mu$ denote the corresponding 
one-dimensional $P$-module. Consider
$\bU=G\times_P (N\oplus \BC_\mu)$ with
projections $\pi: \bU\to G\times_P N$ and
$\kappa: \bU\to G/P$. Then $\pi$ makes $\bU$  
the total space of a line bundle on 
$\bZ$. For simplicity, 
the  sheaf of sections of this bundle is often denoted by $\cl_\bZ(\mu)$
in place of $\cl_\bZ(\BC_\mu)$. Note that  $\cl_\bZ(\mu)^\star=\cl_\bZ(-\mu)$. 
We regard $\BC_\mu$ as the highest weight
space in the $G$-module $V_\mu$. Therefore $\bU$ admits the collapsing into 
$V\oplus V_\mu$.

Since $\bU$ is the total space of a $G$-linearised vector bundle on $G/P$, the 
canonical bundle $\omega_\bU$ is a pull-back of a line bundle on $G/P$.
The top exterior power of the cotangent  space at $e\ast \tilde n\in \bU$ ($e\in G$ is the identity and $\tilde n\in N\oplus \BC_\mu$)
is 
\[
   \wedge^{top}(\g/\p)^*\otimes  \wedge^{top}N^* \otimes (\BC_\mu)^*=
   \wedge^{top}\n \otimes (\wedge^{top}N)^* \otimes (\BC_\mu)^*.
\]
The corresponding character of $P$ is $\gamma-\mu$, where
$\gamma:=|\Delta(\n)|-|\Psi|$.
Therefore  
\[
  \omega_\bU\simeq \kappa^\ast\bigl(\cl_{G/P}(\BC_{\gamma-\mu})\bigr)\simeq
  \pi^\ast\bigl(\cl_\bZ(\gamma-\mu)\bigr).
\]
By Lemma~\ref{lm:odin}, we obtain  $\pi_\ast(\omega_\bU)=
\bigoplus_{n\ge 0} \cl_\bZ(\gamma-\mu)\otimes  \cl_\bZ(n\mu)^\star$ and hence
\[
   H^i(\bU,\omega_\bU)=\bigoplus_{n\ge 0} H^i(\bZ,\cl_\bZ((n{+}1)\mu{-}\gamma)^\star) .
\]
In order to apply Theorem~\ref{thm-GR}, we need sufficient conditions for the
collapsing 
\[
       f_\mu: \bU\to G{\cdot}(N\oplus \BC_\mu)
\]
to be generically finite.  There are two possibilities now.

\textbf{A) \textsf{The collapsing $f:  \bZ\to G{\cdot}N$ is generically finite.}} \\
It is then easily seen that $f_\mu$ is generically finite for any $\mu\in\Xp$. This yields
the following vanishing result:

\begin{prop}    \label{pr:coll-gf}
If $f_0:  \bZ\to G{\cdot}N$ is generically finite and $\gamma=|\Delta(\n)|-|\Psi|$, then 
\[
   H^i(\bZ,\cl_\bZ((n{+}1)\mu{-}\gamma)^\star)=0
\]
for any $\mu\in\mathfrak X^P_+$ and all $n\ge 0$, $i\ge 1$.  
In particular, taking $n=0$ and letting $\nu=\mu-\gamma$, we obtain
\[
 \text{$H^i(\bZ,\cl_\bZ(\nu)^\star)=0$ \ for all $i\ge 1$}
\]
if  $\nu\in \mathfrak X^P$ is such that $\nu \gtrdot |\Psi|-|\Delta(\n)| $.
\end{prop}

\textbf{B)  \textsf{The collapsing $f:  \bZ\to G{\cdot}N$ is \un{not} 
generically finite.}} \\
Here we have to correct the situation, i.e., choose $\mu$ such that $f_\mu$ 
to be generically finite.
 
Looking at the collapsing 
$f_\mu:  G\times_P (N\oplus \BC_\mu)\to G{\cdot}(N\oplus \BC_\mu)$
the other way around, we notice that if $\psi_\mu: G\times_P \BC_\mu\to G{\cdot}\BC_\mu\subset
V_\mu$ is generically finite, then so is $f_\mu$. However, $\psi_\mu$ is generically
finite (in fact, birational) if and only if $\mu\in\mathfrak X^P_+$ is a $P$-regular dominant weight,
i.e., $\mu \gtrdot \rho_P$.  Equivalently, $\mu=\tilde\mu+\rho_P$ for some
$\tilde\mu\in \mathfrak X^P_+$.

This provides a weaker vanishing result that applies 
to arbitrary $P$-submodules. 

\begin{prop}  \label{pr:coll-ngf}
Let $N$ be an arbitrary $P$-submodule.
If $\mu\in \mathfrak X^P_+$ and $\mu\gtrdot \rho_P$, then
\[
   H^i(\bZ,\cl_\bZ((n{+}1)\mu{-}\gamma)^\star)=0
\]
for all $n\ge 0$, $i\ge 1$. In particular, taking $n=0$ and letting 
$\nu=\mu-\gamma$, we obtain 
\[
  \text{$H^i(\bZ,\cl_\bZ(\nu)^\star)=0$ for all $i\ge 1$}
\]
whenever $\nu\in \mathfrak X^P$ and $\nu \gtrdot \rho_P+|\Psi|-|\Delta(\n)| $.
\end{prop} 

Combining Propositions~\ref{pr:coll-gf} and \ref{pr:coll-ngf}, we obtain Theorem~\ref{rez-2}.

\begin{rmk}  The estimate in part B) is not optimal, because  we do not actually
need generic finiteness for $\psi_\mu$. It can happen that both $f$ and $\psi_\mu$
are not generically finite, while $f_\mu$ is. (See e.g. Theorem~\ref{thm:estim-short} below.)
\end{rmk}

\subsection{Proof of Theorem~\ref{rez-1}}
The cohomology groups of $\cl_\bZ(\mu)=\cl_{G\times_P N}(\mu)$ have a natural structure 
of a graded $G$-module by
\[
    H^i(G\times_P N, \cl_{G\times_P N}(\mu))\simeq \bigoplus_{j=0}^\infty
    H^i(G/P, \cl_{G/P}(\mathcal S^j N^*\otimes \BC_\mu)) ,
\]
where $\mathcal S^j N^*$ is the $j$-th symmetric power of the dual of $N$. 
Set $H^i(\mu):=H^i(\bZ, \cl_{\bZ}(\mu)^\star)$.
It is a graded $G$-module with
\[
   (H^i(\mu))_j=H^i(G/P, \cl_{G/P}(\mathcal S^j N\otimes \BC_\mu)^\star).
\]
As $\dim (H^i(\mu))_j < \infty$, the $G$-Hilbert series of 
$H^i(\mu)$ is well-defined:
\[
     \mathcal H_{G}( H^i(\mu); q)=\sum_j \sum_{\lb\in\Xp} \dim \Hom_G(V_\lb, (H^i(\mu))_j)
     e^\lb q^j \in \BZ [\mathfrak X] [[q]] .
\]
We also need the non-graded version of functor $\mathcal H_G$. If $M$ is a finite-dimensional
$G$-module, then 
\[
    \mathcal H_{G}(M)=\sum_{\lb\in\Xp} \dim \Hom_G(V_\lb, M) e^\lb \in \BZ [\mathfrak X].
\]
This extends to virtual $G$-modules by linearity.

Assume for a while that  $P=B$, i.e., $\bZ=G\times_B N$.
By the Borel-Weil-Bott theorem for $G/B$, we have
\[
   H^i(G/B,\cl_{G/B}(\mu)^\star)=
   \begin{cases}  V_\nu^*, & \text{if $\nu=w(\mu+\rho)-\rho \in \Xp$ and $\ell(w)=i$.} \\ 
   0,  & \text{otherwise}. 
   \end{cases}
\]
Using the non-graded functor $\mathcal H_G$, one can also write
\begin{equation}   \label{eq:BWB}
   \mathcal H_G(\sum_i (-1)^i H^i(G/B,\cl_{G/B}(\mu)^\star))=
   \begin{cases}   \esi(w) e^{\nu^*}, &  \text{if $\nu=w(\mu+\rho)-\rho \in \Xp$.}\\
    0,  & \text{otherwise}. 
 \end{cases}
\end{equation}

The following result is well known in case of Lusztig's $q$-analogues, see e.g. 
\cite[Lemma\,6.1]{rkb89}. For convenience of
the reader, we provide a proof of the general statement.

\begin{thm}   \label{thm:svyaz-B}
For any $\mu\in\mathfrak X$, we have 
\[
   \displaystyle\sum_{i}  (-1)^i \mathcal H_G\bigl(H^i(G\times_B N,\cl_{G\times_B N}(\mu)^\star);q\bigr)=\sum_{\lb\in\Xp}\me_{\lb,\Psi}^\mu(q) e^{\lb^*} .
\]
\end{thm}
\begin{proof}
Each finite-dimensional $B$-module $M$ has a 
$B$-filtration such that
the associated graded $B$-module, denoted $\widetilde{M}$, is completely 
reducible.  Then 
\[
   \sum_i (-1)^i H^i(G/B, \cl_{G/B}(M)^\star)=
    \sum_i (-1)^i H^i(G/B, \cl_{G/B}((\widetilde{M})^\star) .
\]
We will apply this to the $B$-modules $\mathcal S^jN\otimes \BC_\mu$, $j=0,1,\dots$.
\begin{multline*}
\sum_{i}  (-1)^i \mathcal H_G(H^i(G\times_B N,\cl_{G\times_B N}(\mu)^\star; q)
\\
=\sum_{j=0}^\infty \mathcal H_G\bigl(\sum_i (-1)^i H^i(G/B, \cl_{G/B}(\mathcal S^j N\otimes \BC_\mu)^\star);q\bigr) 
\\
=\sum_{j=0}^\infty \mathcal H_G\bigl(\sum_i (-1)^i H^i(G/B, \cl_{G/B}(\widetilde{\mathcal S^j N}\otimes \BC_\mu)^\star);q\bigr)
\\
=\sum_{j=0}^\infty \sum_{\nu\vdash \mathcal S^j N}\dim (\mathcal S^j N)^\nu q^j{\cdot} 
\mathcal H_G\bigl( \sum_i (-1)^i H^i(G/B,\cl_{G/B}(\nu+\mu)^\star   \bigr)
\\
=\sum_{\nu\vdash \mathcal S^\bullet N} \eus P_{\Psi,q}(\nu) 
\mathcal H_G\bigl( \sum_i (-1)^i H^i(G/B,\cl_{G/B}(\nu+\mu)^\star   \bigr),
\end{multline*}
where notation $\nu\vdash \mathcal S^j N$ means that $\nu$ is a weight of $\mathcal S^j N$.
By the BWB-theorem, the weight $\nu+\mu$ contributes to the last sum if and only if
$\nu+\mu+\rho$ is regular, i.e., $w(\nu+\mu+\rho)-\rho=\lb\in \Xp$ for a unique $w\in W$.
Therefore,  using Eq.~\eqref{eq:BWB}, we obtain
\begin{multline*}
\sum_{\nu} \eus P_{\Psi,q}(\nu) 
\mathcal H_G\bigl( \sum_i (-1)^i H^i(G/B,\cl_{G/B}(\nu+\mu)^\star   \bigr)=
\\
\sum_{\lb\in\Xp}\sum_{w\in W} \esi(w) \eus P_{\Psi,q}(w^{-1}(\lb+\rho)-\mu-\rho) e^{\lb^*}=\sum_{\lb\in\Xp} \me_{\lb,\Psi}^\mu(q) e^{\lb^*},
\end{multline*}
as required.
\end{proof}

\begin{thm}   \label{thm:svyaz-P}
For any $\mu\in\mathfrak X^P$, we have 
\[
   \displaystyle\sum_{i}  (-1)^i \mathcal H_G\bigl(H^i(G\times_P N,\cl_{G\times_P N}(\mu)^\star);q\bigr)=\sum_{\lb\in\Xp}\me_{\lb,\Psi}^\mu(q) e^{\lb^*} .
\]
\end{thm}
\begin{proof}
Using the Leray spectral sequence associated to the morphism $G/B\to G/P$, one easily proves that, for any $\mu\in \mathfrak X^P$, there is an isomorphism
\[
   H^i(G/B,\cl_{G/B}(\mathcal S^j N\otimes \BC_\mu)^\star)\simeq 
   H^i(G/P,\cl_{G/P}(\mathcal S^j N\otimes \BC_\mu)^\star) .
\]
Thus, the assertion reduces to the previous theorem.
\end{proof}

\begin{cor}   \label{cor:non-neg}
If $\mu\in \mathfrak X^P$ and $H^i(G\times_P N,\cl_{G\times_P N}(\mu)^\star)=0$ for $i\ge 1$, then $\me_{\lb,\Psi}^\mu(q)$ has non-negative coefficients for all $\lb\in \Xp$.
\end{cor}

Now, combining this corollary and Propositions~\ref{pr:coll-gf},\,\ref{pr:coll-ngf},
we obtain Theorem~\ref{rez-1}.

\begin{rmk}     \label{rmk:higher-vanish}
By Theorem~\ref{thm:svyaz-P}, if higher cohomology of $\cl_\bZ(\mu)^\star$ vanishes, then
the polynomial 
$\me_{\lb,\Psi}^\mu(q)$ counts occurrences of $V_\lb^*$ in the graded $G$-module
$H^0(\bZ, \cl_\bZ(\mu)^\star)$.
In particular, $\me_{\lb,\Psi}^\mu(1)$ is the multiplicity of $V_\lb^*$ in 
$H^0(\bZ, \cl_\bZ(\mu)^\star)$.
\end{rmk}

\subsection{}
If we wish to get a generically finite collapsing for a  $B$-stable 
$N\subset V$, then $P$ must be chosen as large as possible. That is, we have to take
$P=\mathsf{Norm}_G(N)$, the normaliser of $N$ in $G$. 
However, even this does not guarantee the generic finiteness.

\begin{ex}    \label{ex:dim-odd}
Let $\ce$ be a $B$-stable subspace of $\ut\subset \g$. 
Actually, $\ce$ is a $B$-stable ideal of $\ut$. Let $P=\mathsf{Norm}_G(\ce)$.
The image of the collapsing $G\times_P\ce\to G{\cdot}\ce$ is the closure of a nilpotent
orbit. Hence $\dim (G{\cdot}\ce)$ is even. However, 
$\dim (G\times_P\ce)$ can be odd. For instance, take $\ce=[\ut,\ut]$. If $G$ is simple and
$G\ne SL_2$, then $\mathsf{Norm}_G([\ut,\ut])=B$. But $\dim (G\times_B [\ut,\ut])$ is 
even if and only if $\rk(G)$ is. It can be shown that the collapsing
$G\times_B [\ut,\ut] \to G{\cdot}[\ut,\ut] $ is generically finite if and only if $\g\in\{
\GR{A}{2n},\,\GR{B}{2n},\,\GR{C}{2n},\,
\GR{E}{6},\,\GR{E}{8},\,\GR{F}{4},\,\GR{G}{2}\}$.
\end{ex}

$B$-stable (or ``\textrm{ad}-nilpotent'') ideals of $\ut$ provide the most natural class of examples of generalised  Kostka-Foulkes polynomials. There is a rich combinatorial theory of these ideals.
In particular, the normalisers of \textrm{ad}-nilpotent ideals has been studied in \cite{ya-norm}.

\begin{ex}
a) For $G=SL_{2n+1}$, consider $\Psi=\{\gamma\in \Delta^+\mid \hot(\gamma)\ge n+1\}$. 
The corresponding \textrm{ad}-nilpotent ideal is $\ut_n=\underbrace{[\dots [}_{n}\ut,\ut],\dots,\ut]$.  By direct calculations, $|\Psi|=\rho$. Therefore 
the normaliser of $\ut_n$ equals $B$ \cite[Theorem\,2.4(ii)]{ya-norm}. 
Next, 
$\dim (G\times_B \ut_n)=2n^2+2n +\genfrac{(}{)}{0pt}{}{n}{2}$  and 
the dense orbit in $G{\cdot}\ut_n$ corresponds to the partition $(2,\dots,2,1)$.
Therefore  
$\dim G{\cdot}\ut_n=2n^2+2n$,  and the collapsing is not generically finite unless
$n=1$. By Theorems~\ref{rez-2}(i) and \ref{rez-1}(i) with $P=B$, we obtain

{\bf --} \ $H^i(G\times_B \ut_n, \cl_{G\times_B \ut_n}(\mu)^\star)=0$ for any $\mu\in\Xp$ and $i\ge 1$;

{\bf --} \ $\me_{\lb,\Psi}^\mu(q)$ has non-negative coefficients for all $\lb, \mu\in\Xp$.

b)  For $G=SL_{2n}$, consider $\Psi=\{\gamma\in \Delta^+\mid \hot(\gamma)\ge n\}$. 
The corresponding \textrm{ad}-nilpotent ideal is $\ut_{n-1}$. 
Since $|\Psi|=\rho+\vp_n$, the normaliser of $\ut_{n-1}$ equals $B$. 
Again, direct calculations show that
$\dim (G\times_B \ut_{n-1})-  \dim G{\cdot}\ut_{n-1}=\genfrac{(}{)}{0pt}{}{n}{2}$.
Here we have

{\bf --} \ $H^i(G\times_B \ut_n, \cl_{G\times_B \ut_{n-1}}(\mu)^\star)=0$ for any $\mu\gtrdot \vp_n$ and $i\ge 1$;

{\bf --} \ $\me_{\lb,\Psi}^\mu(q)$ has non-negative coefficients for all $\lb\in\Xp$ and $\mu\gtrdot \vp_n$.
\end{ex}

\begin{rmk}
For an arbitrary $B$-stable subspace $N\subset V$, the normaliser of $N$ is fully 
determined by
$|\Psi|$. The proof of  \cite[Theorem\,2.4(i),(ii)]{ya-norm} goes thorough verbatim, and it
shows that $|\Psi|$ is dominant and
\[
\left\{  \text{\parbox{160pt}{
the root subspace $\g_{-\ap}$  ($\ap\in\Pi$) belong to 
$\Lie (\mathsf{Norm}_G(N))$} } 
\right\}   \Leftrightarrow  
\left\{  
(\ap, |\Psi|)=0 \right\}.
\]
Equivalently, one can say that  $\mathsf{Norm}_G(N)=\mathsf{Norm}_G(\wedge^{\dim N} N)$,
where  
$\wedge^{\dim N}  N \subset \wedge^{\dim N}  V$.
\end{rmk}

\section{The little adjoint module and short $q$-analogues} 
\label{sect:short}

Let $G$ be a simple algebraic group such that $\Delta$ has
two root lengths. There is a special interesting case in which $\Psi=\Delta^+_s$ is the set 
of short  positive roots.  The subscripts
`s` and `l` will be used to mark objects related to short and long roots, respectively.
For instance, $\Delta_l$ is the set of all long roots, $\Delta^+=\Delta^+_s\sqcup
\Delta^+_l$, and $\Pi_s=\Pi\cap \Delta_s$. 
Let $\bar\theta$ be the short dominant root. The $G$-module $\vbt$ is said to be 
{\it little adjoint}.

\begin{lm}    \label{lm:vesa-little}
The set of nonzero weights of $\vbt$ is $\Delta_s$;
$m_{\bar\theta}^\nu=1$ for  $\nu\in \Delta_s$ and $m_{\bar\theta}^0=\# \Pi_s$.
\end{lm}
The last equality is proved in \cite[Prop.\,2.8]{spin}; the rest is obvious.
It follows that there is a unique $B$-stable subspace
of $\vbt$ whose set of weights is $\Delta^+_s$. Write $\vbt^+$ for this
subspace.  In the rest of the article, we work with 
$\Psi=\Delta^+_s$ and the $B$-stable subspace $N=\vbt^+$.
In place of $\eus P_{\Delta^+_s,q}(\nu)$ and $\me^\mu_{\lb,\Delta^+_s}(q)$, 
we write $\ov{\eus P}_q(\nu)$ and $\ov{\me}_{\lb}^{\mu}(q)$, respectively. 
The polynomials $\ov{\me}_{\lb}^{\mu}(q)$
are said to be {\it short $q$-analogues} (of weight multiplicities).

We have $\Xp\cap \Delta^+_s=\{\bar\theta\}$.
Set $\rho_s=\frac{1}{2}\vert \Delta^+_s\vert$ and $\rho_l=\frac{1}{2}\vert
\Delta^+_l\vert$. It is easily seen that $\rho_s$ (resp. $\rho_l$) is the sum of fundamental
weights corresponding to $\Pi_s$ (resp. $\Pi_l$).
Let $H$ be the connected semisimple subgroup of $G$ that contains $T$ and whose root system is $\Delta_l$. The Weyl group of $H$ is the normal subgroup of $W$  generated
by all "long" reflections. It is denoted by $W_l$. 
Let  $G(\Pi_s)$ (resp. $\g(\Pi_s)$) denote  the 
simple subgroup of $G$ (subalgebra of $\g$) whose set of simple roots is $\Pi_s$.
Then $\rk\,\g(\Pi_s)=\# \Pi_s$ and $B\cap G(\Pi_s)=:B(\Pi_s)$ is a Borel subgroup of 
$G(\Pi_s)$. Clearly,  $G(\Pi_s){\cdot}T=:L$ is a standard Levi subgroup of $G$ and
 $G(\Pi_s)= (L,L)$.

The collapsing $G\times_B \vbt^+\to 
G{\cdot}\vbt^+$ is not generically finite, and
Theorem~\ref{rez-1}(i) (with $\rho_P=\rho$, $\Delta(\n)=\Delta^+$, and $|\Delta^+_s|=
2\rho_s$) yields the bound $\mu\gtrdot 2\rho_s-\rho=\rho_s-\rho_l$ for 
$\ov{\me}_{\lb}^{\mu}(q)$.
However, in this case there is a better bound, and our first goal is to obtain it.
To this end, we need some further properties of  little adjoint modules. 

The weight structure of $\vbt$ shows that  
$\vbt\vert_{G(\Pi_s)}$ contains the adjoint representation of $G(\Pi_s)$. 
To distinguish the Lie algebra $\g(\Pi_s)$ sitting in $\g$ and the adjoint representation of $G(\Pi_s)$ sitting in $\vbt$, the latter will be denoted by $\widehat{\g}(\Pi_s)$.
That is, 
\[
    \vbt\vert_{G(\Pi_s)}=\widehat{\g}(\Pi_s)\oplus  R , 
\]
where $R$ is {\sl the\/} complementary $G(\Pi_s)$-submodule. The above decomposition
is $L$-stable and hence $T$-stable. We have $R^T=0$ and the weights of $R$ are those 
short roots that are not $\BZ$-linear combinations of short simple roots.
Furthermore, 
$\vbt^+=\widehat{\g}(\Pi_s)^+\oplus  R^+$, where $R^+\subset R$ and 
$\g(\Pi_s)^+=\g(\Pi_s)\cap\ut$  is a maximal nilpotent subalgebra of $\g(\Pi_s)$.

\begin{thm}     \label{thm:estim-short}
If $\mu\gtrdot \rho_l$, then the collapsing
$f_\mu^{(s)}:  G\times_B (\vbt^+\oplus \BC_\mu) \to  G{\cdot}(\vbt^+\oplus \BC_\mu)$ is birational.
\end{thm}
\begin{proof}  
Recall that $\BC_\mu$ is the line of $B$-highest weight vectors in $V_\mu$.
Obviously, $f_\mu^{(s)}$ is  birational if and only if the following property holds: for a generic point $(v,v_\mu)\in \vbt^+\oplus \BC_\mu$,
if $g{\cdot}(v,v_\mu)\in \vbt^+\oplus \BC_\mu$ ($g\in G$), then $g\in B$.
Let  $\tilde P$  denote the standard parabolic subgroup of $G$ whose Levi subgroup is $L$. 
If $\mu\gtrdot \rho_l$, then
the normaliser in $G$ of the line $\langle v_\mu\rangle$ is contained in $\tilde P$.
Consequently, if $g{\cdot}(v,v_\mu)\in \vbt^+\oplus \BC_\mu$, then $g\in \tilde P$.

Take $v=v'+ r\in \vbt^+$ ($r\in R$) such that $v'$ is a regular nilpotent element
of $\widehat{\g}(\Pi_s)^+$.  Write $g=g_1g_2\in\tilde P$, where $g_1\in G(\Pi_s)$ and $g_2$ lies in the radical
of  $\tilde P$, $\mathsf{rad}(\tilde P)$. It is easily seen that $\mathsf{rad}(\tilde P)$ preserves
$R^+$ and acts trivially in $\vbt^+/R^+$. Therefore $g_2$
does not change the $\widehat{\g}(\Pi_s)$-component of
$v$, i.e., $g_2{\cdot}v=v'+r'$ ($r'\in R^+$).  Hence $g{\cdot}v=g_1{\cdot}v' + g_1{\cdot}r'$,
and $g_1{\cdot}v'\in \widehat{\g}(\Pi_s)^+$ is still a regular nilpotent element of 
$\widehat{\g}(\Pi_s)$.
But the latter is only possible if $g_1\in B(\Pi_s)$ and hence $g\in B$.
\end{proof}

\begin{cor}   \label{cor:rho-l}
If $\nu+\rho_l\in\Xp$, then 
\begin{itemize}
\item[\sf (i)] \ 
$H^i(G\times_B V^+_{\bar\theta}, \cl_{G\times_B V^+_{\bar\theta}}(\nu)^\star)=0$
for $i\ge 1$;   
\item[\sf (ii)] \ $\ov{\me}_{\lb}^\nu(q)$ has non-negative coefficients for all
$\lb\in\Xp$.
\end{itemize}
\end{cor}
\begin{proof}
(i) Set $\bU=G\times_B (\vbt^+\oplus \BC_\mu)$ and $\bZ=G\times_B V^+_{\bar\theta}$. Then $\omega_\bU=\cl_\bU(\gamma-\mu)$, where $\gamma=|\Delta^+|- |\Delta^+_s|=2\rho_l$. By Theorems~\ref{thm-GR} and \ref{thm:estim-short}, 
$H^i(\bU,\omega_\bU)=0$ for $i\ge 1$ whenever $\mu\gtrdot \rho_l$. 
Hence $H^i(\bZ,\cl_\bZ((n+1)\mu-\gamma)^\star)=0$, see Section~\ref{sect:cohomol}.
In particular, $H^i(\bZ,\cl_\bZ(\nu)^\star)=0$, where $\nu=\mu-\gamma$.
It remains to observe that $\nu\gtrdot -\rho_l$.

(ii)  This follows from (i) and Theorem~\ref{thm:svyaz-B}.
\end{proof}

\begin{rmk}
The  proof of  Corollary~\ref{cor:rho-l}(i) uses 
(a version of) the Grauert--Riemenschneider theorem.
However, for $\nu=0$ (at least) 
one can adapt Hesselink's proof of \cite[Theorem\,B]{wim1}, which does 
not refer to Grauert--Riemenschneider and goes through for any algebraically closed field 
$\bbk$ of characteristic zero.
Using this, one can prove the following:
Let $\tilde N$ be any $B$-stable subspace of $\vbt$ such that 
$\tilde N\supset  \vbt^+$. Then
$H^i(G\times_B \tilde N, \co_{G\times_B \tilde N})=0$ for $i\ge 1$.
\end{rmk}

\noindent
Let us describe a semi-direct product structure of $W$, which plays an important role
below. Consider two subgroups of $W$:

\textbullet \ \ $W_l$ is generated by \un{all} ``long" reflections in $W$. 
It is a normal subgroup of $W$.

\textbullet \ \ $W(\Pi_s)$ is generated by all \un{sim}p\un{le} ``short" reflections, i.e., by 
$s_{\ap}$ with $\ap\in \Pi_s$.

\begin{lm}   \label{lm:semi-dir}  
\begin{itemize}
\item[\sf (i)] \  $W$ is a semi-direct product of $W_l$ and $W(\Pi_s)$: 
$W\simeq W(\Pi_s)\ltimes W_l$.
\item[\sf (ii)] \ $W(\Pi_s)=\{w\in W \mid w(\Delta^+_l)\subset \Delta^+_l\}$.
\end{itemize}
\end{lm}
\begin{proof}
(i) Since $W_l$ is a normal subgroup of $W$ and $W_l\cap W(\Pi_s)=\{1\}$, 
it suffices to prove that the natural mapping 
$W(\Pi_s)\times W_l\to W$ is onto. We argue by induction on the length of $w\in W$.
Suppose $w\not\in W(\Pi_s)$ and
$w=w_1s_\beta w_2\in W$, $\beta\in\Pi_l$, is a reduced decomposition.
Then $w=w_1w_2s_{\beta'}$, where $\beta'=w_2(\beta)\in\Delta_l$,
and $\ell(w_1w_2)<\ell(w)$. 
Thus, all long simple reflections occurring in an expression for $w$
can eventually be moved up to the right.

(ii)  Since 
$s_\ap(\Delta^+_l)\subset \Delta^+_l$ for $\ap\in\Pi_s$, 
$W(\Pi_s)\subset\{w\in W \mid w(\Delta^+_l)\subset \Delta^+_l\}$.
On the other hand, if $w(\Delta^+_l)\subset \Delta^+_l$ and
$w=w's_\ap$ is a reduced decomposition, then the equality $\mathsf{N}(w)=s_\ap(\mathsf{N}(w'))\cup\{\ap\}$
shows that $\ap$ is necessarily short, so that
we can argue by induction on $\ell(w)$.
\end{proof}

Recall that the {\it null-cone\/} of a $G$-module $V$, $\fN(V)$, is the  
zero set of all homogeneous $G$-invariant polynomials of positive degree.
Next proposition summarises  invariant-theoretic properties of $\vbt$ and
$\fN(\vbt)$ required below, which are of independent interest.
All the assertions can easily be verified using the classification, but our intention 
is to present a conceptual proof.

\begin{prop}   \label{nullcone-vbt}
a) \ $\fN(\vbt)=G{\cdot}\vbt^+$. Hence it is irreducible; 

b) \  The restriction homomorphisms $\BC[\vbt] \to \BC[\widehat{\g}(\Pi_s)] \to \BC[\vbt^0]$
induce the isomorphisms
$\BC[\vbt]^G\isom \BC[\widehat{\g}(\Pi_s)]^{G(\Pi_s)}
\isom \BC[\vbt^0]^{W(\Pi_s)}$, and  $\BC[\vbt]^G$ is a polynomial algebra.

c) \  $\fN(\vbt)$ is a reduced normal complete intersection of codimension $\#(\Pi_s)$.
\end{prop}
\begin{proof}[Outline of the proof] We refer to \cite{vipo} for invariant-theoretic results 
mentioned below. 

a) This follows from the Hilbert-Mumford criterion and the fact any maximal subset of weights
of $\vbt$, lying in an open half-space, is $W$-conjugate to $\Delta^+_s$.

b) The weight structure of $\vbt$ shows that $\vbt^0=\vbt^H$. If $v\in \vbt^0$
is generic, then $\g{\cdot}v+\vbt^0=\vbt$. Therefore $G{\cdot}\vbt^0$ is dense in $\vbt$
and a generic stabiliser (= {\it stabiliser in general
position}) for $G{:}\vbt$ contains $H$. Actually, it is not hard to prove that
$H$ is a generic stabiliser for $G{:}\vbt$. By the Luna-Ri\-chard\-son theorem, we then have
$\BC[\vbt]^G\isom \BC[\vbt^H]^{N_G(H)/H}$, and it is easily seen that  
$N_G(H)/H\simeq W/W_l\simeq W(\Pi_s)$. Furthermore, the $W(\Pi_s)$-action on
$\vbt^0$ is nothing but the standard reflection representation on the Cartan subalgebra
of $\g(\Pi_s)$.

c)  Let $f_1,\dots,f_m$ be basic invariants in $\BC[\vbt]^G\simeq 
\BC[\widehat{\g}(\Pi_s)]^{G(\Pi_s)}$, $m=\#(\Pi_s)$.
Let $e\in \widehat{\g}(\Pi_s)\subset \vbt$ be regular nilpotent. Then
the differentials of the $f_i$'s are linearly independent at 
$e\in \fN(\vbt)$ \cite{ko63}. Hence the ideal of $\fN(\vbt)$ is $(f_1,\dots, f_m)$ and
$\fN(\vbt)$ is a reduced complete intersection 
(cf. \cite[Lemma\,4]{ko63}).
Finally, $\fN(\vbt)$ contains a dense $G$-orbit whose complement is of codimension
$\ge 2$. This yields the normality.
\end{proof}

\noindent
Our ultimate goal is to get a complete characterisation of weights 
$\mu\in\mathfrak X$ such that $\ov{\me}_\lb^\mu(q)$
has nonnegative coefficients for any $\lb\in\Xp$. 
To this end, we exploit a different approach that does not
use vanishing theorems of Section~\ref{sect:cohomol}.

A key observation is that short $q$-analogues obey certain symmetries with respect to the simple reflections
$s_\ap\in W$,  $\ap\in \Pi_l$.  Clearly,  $s_\ap(\Delta^+_s)=\Delta^+_s$. Therefore
       $\ov{\eus P}_q(\nu)=\ov{\eus P}_q(s_\ap \nu)$.
Using this, we compute
\begin{multline}   \label{eq:cert_symm}
\ov{\me}_{\lb}^\mu(q)=\sum_{w\in W} \esi(w)  \ov{\eus P}_q(w(\lb+\rho)-(\mu+\rho)) 
\\
 =\sum_{w\in W} \esi(w)  \ov{\eus P}_q(s_\ap w(\lb+\rho)-s_\ap(\mu+\rho))
 =-\sum_{w\in W} \esi(w)  \ov{\eus P}_q(w(\lb+\rho)-s_\ap\mu-s_\ap\rho)
\\
= -\sum_{w\in W} \esi(w)  \ov{\eus P}_q(w(\lb+\rho)-(s_\ap\mu-\ap+\rho))
=- \ov{\me}_{\lb}^{s_\ap(\mu+\ap)}(q) .
\end{multline}

\noindent
The {\it shifted action\/} of $W_l$ on $\mathfrak X$ is defined by
\[
      w\odot \gamma=w(\gamma+\rho_l)-\rho_l .
\]
For $\ap\in\Pi_l$,  one easily recognise $s_\ap(\mu+\ap)$  as $s_\ap\odot \mu$ and
hence Eq.~\eqref{eq:cert_symm} can be written as
    $\ov{\me}_\lb^{s_\ap\odot \mu}(q)=- \ov{\me}_\lb^{\mu}(q)$.
This readily implies the equality
\begin{equation}   \label{eq:shifted_act}
   \ov{\me}_\lb^{w\odot \mu}(q)=\esi(w)\,\ov{\me}_\lb^{\mu}(q)
\end{equation}
for any $w\in W_l$. 
Note that for $w\in W_l$, the length $\ell(w)$  depends on the choice of ambient group, $W$ or $W_l$,
but the parity $\esi(w)$ does not! (This is because $\esi(w)=\det (w)$ for the reflection
representation of $W$ in $\mathfrak X\otimes_\BZ\BQ$.)

Let $\mathfrak X_{+,H}$ denote the monoid of $H$-dominant weights
with respect to $\Delta^+_l$. From \eqref{eq:shifted_act}, we immediately deduce that 

\begin{itemize}
\item 
it suffices to know \/$\ov{\me}_\lb^{\mu}(q)$ for $\mu\in \mathfrak X_{+,H}- \rho_l$.
\item 
if such a $\mu$ is not $H$-dominant, then it lies on a wall of the shifted
dominant Weyl chamber for $H$, and hence $\ov{\me}_\lb^{\mu}(q)\equiv 0$.
\item  Thus, the problem is reduced
to studying polynomials $\ov{\me}_\lb^{\mu}(q)$ for $\mu\in \mathfrak X_{+,H}$.
\end{itemize}

\noindent
Short $q$-analogues enjoy several good interpretations at $q=1$.
Write $\ov{m}_\lb^\nu$ in place of $\ov{\me}_\lb^\nu(1)$.

\begin{enumerate}
\item
As already observed in Remark~\ref{rmk:higher-vanish}, if higher cohomology 
of $\cl_\bZ(\nu)^\star$ vanish, then 
$\ov{m}_\lb^\nu$ is the multiplicity of $V_\lb^*$ in
$H^0(\bZ, \cl_\bZ(\nu)^\star)$.
\item If $\nu\in\mathfrak X_{+,H}$ and $V^{(H)}_\nu$ is a simple $H$-module with 
highest weight $\nu$, then $\ov{m}_\lb^\nu$ is the multiplicity
of $V^{(H)}_\nu$ in $V_\lb\vert_H$, denoted $\mathsf{mult}(V^{(H)}_\nu, V_\lb\vert_H)$,
see  \cite[Lemma\,3.1]{heck}.
\end{enumerate}

\noindent
(Our $\ov{m}_\lb^\nu$ is $m_\lb^{G,H}(\nu)$ in the notation of \cite{heck}.
In fact, Heckman works in a general situation, where  $H\subset G$ is an arbitrary connected reductive group.)
Furthermore, the numbers $\ov{m}_\lb^\nu$ are naturally defined for all $\lb,\nu\in\mathfrak X$
and they satisfy the relation 
\begin{equation}    \label{eq:short-symm-heck}
      \ov{m}_{w(\lb+\rho)-\rho}^{\bar{w}(\nu+\rho_l)-\rho_l}=\esi(w)\esi(\bar w)\ov{m}_\lb^\nu ,
      \quad w\in W,\ \bar w\in W_l.
\end{equation}
(See Equation~(3.7) in \cite{heck}.) 
The semi-direct product structure of $W$ provides an extra symmetry to this picture 
that is absent in the general setting of \cite{heck}.
Namely, if $\nu$ is $H$-dominant, then so is $w\nu$ for any $w\in W(\Pi_s)$.
Using this one easily proves that $\ov{m}_\lb^\nu=\ov{m}_\lb^{w\nu}$ 
for all $\lb \in \Xp$ and $w\in W(\Pi_s)$.

Recall that $\{\mu^+\}=W\mu\cap\Xp$.
Let $w_\mu$ denote the unique element of minimal length such that $w_\mu(\mu)=\mu^+$.

\begin{lm}    \label{lm:pro-w-mu}
If $\mu \in \mathfrak X_{+,H}$, then $w_\mu\in W(\Pi_s)$ and hence 
$\mu^+ -\mu$ is a nonnegative $\BZ$-linear combination of short simple roots.
\end{lm}\begin{proof}
It is known that $\mathsf{N}(w_\mu)=\{\gamma\in \Delta^+\mid (\gamma,\mu)<0\}$, see \cite[Prop.\,2(i)]{br97}.
Since $\mu$ is $H$-dominant, $\mathsf{N}(w_\mu)\subset \Delta^+_s$, and we conclude by 
Lemma~\ref{lm:semi-dir}(ii).
\end{proof}

\begin{prop}    \label{pr:neg-coef}
Let $\mu\in\mathfrak X_{+,H}$.
  
1) Suppose that there is $\nu\in \Xp$ such that 
$\mu\curle \nu\prec \mu^+$. 
Then $\ov{\me}_\nu^\mu(q) \ne 0$
and $\ov{m}_\nu^\mu=0$. In particular, $\ov{\me}_\nu^\mu(q)$ has both positive and negative coefficients.

2) If $V_{\mu^+}^*$ occurs in 
$H^{0}(G/B,\cl_{G/B}(\widetilde{\cs^j(\vbt^+)}\otimes\BC_\mu)^\star)$, then 
$j\ge \hot(\mu^+-\mu)$. Furthermore, for $j=\hot(\mu^+-\mu)$, $H^0(\dots)$ contains a 
unique copy of $V_{\mu^+}^*$.
\end{prop}
\begin{proof}
1) \ Since $w_\mu\in W(\Pi_s)$, we have $\ov{m}_\nu^\mu=
\ov{m}_\nu^{\mu^+}$, and the latter equals zero, because $\nu\prec \mu^+$.
(Obviously, the $H$-module with highest weight $\mu^+$ cannot
occur in $V_\nu\vert_H$.)

Since $\mu\curle \nu\prec \mu^+$ and $\mu^+ -\mu$ is a nonnegative 
$\BZ$-linear combination of short simple roots,  the latter holds for $\nu-\mu$ as well.
Set $a=\hot(\nu-\mu)$.
By definition, 
\[
\ov{\me}_\nu^\mu(q)=\sum_{w\in W}\esi(w) \ov{\eus P}_q(w(\nu+\rho)- (\mu+\rho)) .
\]
As $\nu-\mu \in \mathsf{Span}(\Pi_s)$,  the summand 
$\ov{\eus P}_q(w(\nu+\rho)- (\mu+\rho))$ can be nonzero only if $w\in W(\Pi_s)$.
For $w=1$, we have $\ov{\eus P}_q(\nu-\mu)=q^a +\text{ (lower terms)}$.
If $w\ne 1$, then $\deg \ov{\eus P}_q(w(\nu+\rho)- (\mu+\rho)) < a$. Hence
the highest term of $\ov{\me}_\nu^\mu(q)$ is $q^a$, and we are done.

2) \ This readily follows from the BWB-theorem and Lemma~\ref{lm:pro-w-mu}.
\end{proof}

\noindent
Our main result on non-negativity for short $q$-analogues is a converse to the first claim of 
the previous proposition. For the proof of the main theorem, we need a technical lemma.

\begin{lm}    \label{lm:technical}
{\sf 1)} \  Suppose that $V_\nu^*$ occurs in 
$H^i(G/B, \cl_{G/B}(\wedge^j(\vbt/V^+_{\bar\theta})\otimes \BC_\mu)^\star)$.
Then $\nu\curle \mu^+$.
{\sf 2)} \ (For $\nu=\mu^+$.) If\/ $V_{\mu^+}^*$ occurs in 
$H^i(G/B, \cl_{G/B}(\wedge^j\widetilde{(\vbt/V^+_{\bar\theta}})\otimes \BC_\mu)^\star)$, 
then $j\ge i\ge \ell(w_\mu)$.
\end{lm}
\begin{proof}  Set $M_j=\wedge^j(\vbt/V^+_{\bar\theta})\otimes \BC_\mu$. 

1) If  $V_\nu^*$
occurs in $H^i(G/B, \cl_{G/B}(M_j)^\star)$, then it also occurs in 
$H^i(G/B, \cl_{G/B}(\tilde M_j)^\star)$. By the BWB-theorem, there is then a weight $\gamma$
of $M_j$ and $w\in W$  such that $\ell(w)=i$ and $w(\gamma+\rho)-\rho=\nu$. 
All weights of $M_j$ are of the form
$\mu-|A|$ for some $A\subset \Delta^+_s$, where $\#(A)\le j$.
Hence
$w(\mu+\rho-|A|)=\rho+\nu$. Clearly, $w(\rho-|A|)=\rho-|C|$ for some $C\subset \Delta^+_s$ depending on $w$ and $A$. Thus, $w(\mu+\rho-|A|)\curle w(\mu)+\rho$ and
$\nu\curle w(\mu)\curle \mu^+$.

2) If $V_{\mu^+}^*$ occurs in $H^i(G/B, \cl_{G/B}(\tilde M_j)^\star)$, then, by the 
first part of
the proof, we must have $w(\mu+\rho-|A|)=\rho+\mu^+$, where $A\subset \Delta^+_s$
and $\ell(w)=i$. Hence   $w(\mu)=\mu^+$ and $w(\rho-|A|)=\rho$. 
Therefore $A=\mathsf{N}(w)$ and $i=\ell(w)=\#(A)\ge \ell(w_\mu)$.
Since $\#(A)\le j$ as well, we are done.
\end{proof}

The following is the main result of this section.

\begin{thm}    \label{thm:main}
For $\mu\in\mathfrak X_{+,H}$,  the following conditions are equivalent:
\begin{itemize}
\item[\sf (i)] \  $H^i(G\times_B V^+_{\bar\theta}, \cl_{G\times_B V^+_{\bar\theta}}(\mu)^\star)=0$ for all $i\ge 1$;
\item[\sf (ii)] \  $\ov{\me}_\lb^\mu(q)$ has nonnegative coefficients for any $\lb\in\Xp$;
\item[\sf (iii)] \  If $\mu\curle \nu\curle \mu^+$ for $\nu\in\Xp$, then $\nu=\mu^+$;
\item[\sf (iv)] \  $(\mu,\ap^\vee)\ge -1$ for all $\ap\in\Delta^+_s$.
\end{itemize}
\end{thm}
\begin{proof}
By Corollary~\ref{cor:non-neg}, (i) implies (ii); and Proposition~\ref{pr:neg-coef}
shows that (ii) implies (iii). Since $\nu$ is already assumed to be 
$H$-dominant, (iii) and (iv) are equivalent in view of \cite[Prop.\,2(iii)]{br97}.

It remains to prove the implication (iii) $\Rightarrow$ (i). Our argument is an adaptation
of Broer's proof of \cite[Theorem\,2.4]{br93}. We construct  a similar Koszul complex and consider its spectral sequence of hypercohomology.
 
The pull-back vector bundle $G\times_B (\vbt\oplus (\vbt/\vbt^+))$ on $\bX:=G\times_B \vbt$
has the global $G$-equivariant section  $g\ast v \mapsto g\ast(v, \bar v)$ whose 
scheme of zeros is exactly $\bZ=G\times_B \vbt^+$.
Here $\bar v$ is the image of $v\in \vbt$ in $\vbt/\vbt^+$.
Let $\iota: \bZ\to \bX$ denote the inclusion.
The dual of this section gives rise to a locally free Koszul resolution of $\co_\bZ$ regarded as
$\co_\bX$-module:
\[
   \dots \to \cf^{-1}\to \cf^0 \to \iota_*\co_\bZ \to 0 
\] 
with $\cf^{-j}=\cl_\bX (\wedge^j (\vbt/\vbt^+)^\star[-j]$.
Here the brackets `$[-j]$' denote the degree shift of a graded module.
(That is, if $\mathcal M=\oplus \mathcal M_i$, then $\mathcal M[r]_i=\mathcal M_{r+i}$.)
Therefore the generators of the locally free $\co_\bX$-module $\cf^{-j}$ have degree $j$.
Tensoring this complex with the invertible sheaf $\cl_\bX(\BC_\mu)^\star=
\cl_\bX(\mu)^\star$, we get a locally free resolution of graded $\co_\bX$-modules
\begin{equation}  \label{eq:Koszul}
    \cf(\mu)^\bullet \to \iota_*\cl_\bZ(\mu)^\star\to 0 ,
\end{equation}
where $\cf(\mu)^{-j}=\cl_\bX (\wedge^j (\vbt/\vbt^+)\otimes \BC_\mu)^\star[-j]$. Since
$\bX\simeq G/B\times \vbt$, we have the isomorphism
\[
   H^i(\bX, \cl_\bX (\wedge^j (\vbt/\vbt^+)\otimes \BC_\mu)^\star)\simeq
   \BC[\vbt]\otimes H^i(G/B,\cl_{G/B}(\wedge^j (\vbt/\vbt^+)\otimes \BC_\mu)^\star) 
\]
of graded $\BC[\vbt]$-modules.
For the spectral sequence of hypercohomology associated to the Koszul complex
\eqref{eq:Koszul}, we have 
\[
''E^{kl}_2=H^k(\bX, \eus H^l(\cf(\mu)^\bullet))=
\begin{cases} 
H^k(\bX, \iota_*\cl_\bZ(\mu)^\star)=H^k(\bZ, \cl_\bZ(\mu)^\star),  &  \text{ if  $l=0$}; \\
                                                           0,  &  \text{ if  $l\ne 0$}.   \end{cases}
\]                                                           
and
\[
'E^{kl}_1=H^l(\bX, \cf(\mu)^{k})=
  \BC[\vbt][k]\otimes H^l(G/B, \cl_{G/B}(\wedge^{-k}(\vbt/\vbt^+)\otimes \BC_\mu)^\star) .
\] 
(See \cite[5.7]{weib} for basic facts on hypercohomology.)
It follows that there is a spectral sequence of graded $\BC[\vbt]$-modules
\begin{equation}  \label{eq:spect-seq}
    'E^{-j,i}_1=\BC[\vbt][-j]\otimes H^i(G/B, \cl_{G/B}(\wedge^{j}(\vbt/\vbt^+)\otimes \BC_\mu)^\star) \Rightarrow H^{i-j}(\bZ,\cl_\bZ(\mu)^\star) .
\end{equation}
Let $i-j$ be maximal with $H^i(G/B,\cl_{G/B}(\wedge^j (\vbt/\vbt^+)\otimes \BC_\mu)^\star)\ne 0$.  If $V_\nu^*$ occurs in this cohomology group, then $\nu\curle \mu^+$, by Lemma~\ref{lm:technical}(1).  A basis for $V_\nu^*$ corresponds to some free generators of 
$\BC[\vbt]$-module $'E^{-j,i}_1$ of degree $j$.  Since $i-j$ is maximal, these generators are
in the kernel of $d^{-j,i}_1$. But they are not in the image of  $d^{-j-1,i}_1$, as all elements
of $'E^{-j-1,i}_1$ are of degree $>j$. Hence these generators correspond to nonzero generators of $'E^{-j,i}_2$. Likewise, their images in $'E^{-j,i}_k$ do not vanish.
In view of convergence of the above spectral sequence, this implies that
the multiplicity of $V_\nu^*$ in  $H^{i-j}(\bZ,\cl_\bZ(\mu)^\star)$ is at least one.
It follows that, for some $m\in\BN$, the multiplicity of $V_\nu^*$ in 
$H^{i-j}(G/B,\cl_{G/B}(\widetilde{\cs^m(\vbt^+)}\otimes\BC_\mu)^\star)$ is also at least one.
Any weight of $\cs^m(\vbt^+)\otimes \BC_\mu$ is of the form $\mu+\gamma$ with 
$\gamma\curge 0$. Hence $\nu+\rho=w(\mu+\gamma+\rho)$ for some $w\in W$ with
$\ell(w)=i-j$. Consequently, $\nu\curge \mu$ and altogether $\mu\curle\nu\curle\mu^+$. 
Hence $\nu=\mu^+$. Now, Lemma~\ref{lm:technical}(2) yields $i=j\ge \ell(w_\mu)$. 
In particular,  condition (i) holds.
\end{proof}

The following is an analogue of \cite[Prop.\,2.6]{br93}. 

\begin{prop}
Suppose that $\mu\in\mathfrak X_{+,H}$ satisfies vanishing conditions of Theorem~\ref{thm:main}. Then the graded $\BC[V_{\bar\theta}]$-module 
$H^0(\bZ, \cl_\bZ(\mu)^\star)$ is generated by the unique copy of $V_{\mu^+}^*$ sitting
in degree $\hot(\mu^+-\mu)$.
\end{prop}
\begin{proof}
Eq.~\eqref{eq:spect-seq} an the last part of the proof of Theorem~\ref{thm:main} shows that 
\begin{itemize}
\item The generators of the $\BC[V_{\bar\theta}]$-module $H^0(\bZ, \cl_\bZ(\mu)^\star)$
arise from $G$-modules sitting in 
$H^i(G/B,\cl_{G/B}(\wedge^i (\vbt/\vbt^+)\otimes \BC_\mu)^\star)$, 
with $i\ge \ell(w_\mu)$;
\item  $H^i(G/B,\cl_{G/B}(\wedge^i (\vbt/\vbt^+)\otimes \BC_\mu)^\star)$ only contains
$G$-modules of type $V_{\mu^+}^*$. 
\end{itemize}
It follows that the degree of generators of $H^0(\bZ, \cl_\bZ(\mu)^\star)$ is at least
$\ell(w_\mu)$. On the other hand, if 
$H^{0}(G/B,\cl_{G/B}(\widetilde{\cs^j(\vbt^+)}\otimes\BC_\mu)^\star)$  contains 
a $G$-submodule of type $V_{\mu^+}^*$, then $j\le \hot(\mu^+-\mu)$ by 
Proposition~\ref{pr:neg-coef}(2). Therefore, there cannot be generators of degree larger
than $\hot(\mu^+-\mu)$. It only remains to prove that if $\mu\in\mathfrak X_{+,H}$ satisfies the vanishing condition, then $\ell(w_\mu)=\hot(\mu^+-\mu)$. Clearly, $\ell(w_\mu)\le  \hot(\mu^+-\mu)$. Assume the inequality is strict. Then there is  a $w\in W$ and a simple reflection 
$s_i$ such that $\mu\curle w(\mu) \prec s_iw(\mu)\curle \mu^+$ and 
$s_iw(\mu)=w(\mu)+k\ap_i$ with $k\ge 2$. Then $\nu:=w(\mu)+\ap_i$ belongs to the convex hull of $w(\mu)$ and $s_iw(\mu)$; hence $\mu\prec \nu^+\prec \mu^+$, which
contradicts the vanishing condition.
\end{proof}

Finally, we mention that above two interpretations of numbers $\ov{m}_\lb^\nu$ and
Theorem~\ref{thm:main} lead to an interesting equality.  

\begin{prop}
If $\nu\in \mathfrak X_{+,H}$  and $(\nu,\ap^\vee)\ge -1$ for all $\ap\in\Delta^+_s$, then 
$H^0(G/H, \cl_{G/H}(V^{(H)*}_\nu))$ and 
$H^0(G\times_B V^+_{\bar\theta}, \cl_{G\times_B V^+_{\bar\theta}}(\nu)^\star)$ 
are isomorphic $G$-modules. In particular, for $\nu=0$, we obtain 
$\BC[G/H] \simeq  \BC[G\times_B V^+_{\bar\theta}]$ as $G$-modules.
\end{prop}\begin{proof}
By Frobenius reciprocity, 
\[
   \mathsf{mult}( V_\lb^\ast,  H^0(G/H, \cl_{G/H}((V^{(H)*}_\nu))=
   \mathsf{mult}(V^{(H)}_\nu, V_\lb\vert_H).
\]
Hence the multiplicity of $V_\lb^\ast$ in both spaces $H^0(..)$ under consideration
is equal to $\ov{m}_\lb^\nu$.
\end{proof}

\section{Short Hall-Littlewood polynomials } 
\label{sect:HL}

\noindent
In this section, we define ''short" analogues of Hall-Littlewood polynomials
and establish their basic properties.
Recall that $\Delta$ is a reduced irreducible root system, and
$\Delta^+=\Delta^+_s\sqcup \Delta^+_l$, 
$\Pi=\Pi_s\sqcup \Pi_l$, etc. It is convenient to assume that
in the simply-laced case all roots are short and $\Pi_l=\varnothing$. Then the following  can be regarded as 
a generalisation of Gupta's theory \cite{rkg87, rkg87-A}.

The character ring $\boldsymbol{\Lambda}$ of finite-dimensional
representations of $G$ is identified with $\BZ[\mathfrak X]^W$.
For $\lb\in\Xp$, let 
$\chi_\lb$ denote the character of $V_\lb$,
i.e., $\chi_\lb=\mathsf{ch}(V_\lb)=\sum_\mu m_\lb^\mu e^\mu$. By Weyl's character formula,
$\chi_\lb=  J(e^{\lb+\rho})/ J(e^\rho)$, where $J=\sum_{w\in W} \esi(w)w$
is the skew-symmetrisation operator. Weyl's denominator formula  says that
$J(e^\rho)=e^\rho \prod_{\ap> 0}(1-e^{-\ap})$.
The usual scalar product $\langle\ ,\ \rangle$ on $\boldsymbol{\Lambda}=\BZ[\mathfrak X]^W$ is given by 
$\langle\chi_\lb ,\chi_\nu \rangle=\delta_{\lb,\nu}$.

The  projection $j:  \BZ[\mathfrak X] \to \BZ[\mathfrak X]^W$ is given by
$j(f):= J(f)/J(e^\rho)$.

\noindent 
Set $t_\lb^{(\Pi_s)}(q)=\sum q^{\ell(w)}$, where the summation is over $w\in W(\Pi_s)_\lb$, the 
stabiliser of $\lb$ in $W(\Pi_s)$. 

We will work in the $q$-extended character ring $\boldsymbol{\Lambda}[[q]]$ or its subring
$\boldsymbol{\Lambda}[q]$ and
agree to extend our operators and form $q$-linearly.
We first put
\[
\tilde\bD^{(s)}_q=\frac{ e^\rho}{\prod_{\ap\in\Delta^+_s}(1-qe^\ap)}, 
\quad \bD^{(s)}_q={e^\rho}{\prod_{\ap\in\Delta^+_s}(1-qe^{-\ap})}.
\]
For $\lb,\mu\in\Xp$, define :
\[
  \ov{E}_\mu(q)=j(e^\mu{\cdot}\tilde\bD^{(s)}_q), \quad 
  \ov{P}_\lb(q)=\frac{1}{t_\lb^{(\Pi_s)}(q)}j(e^\lb{\cdot}\bD^{(s)}_q) .
\]
Clearly,  $\ov{E}_\mu(q)\in \boldsymbol{\Lambda}[[q]]$ and $t_\lb^{(\Pi_s)}(q){\cdot}\ov{P}_\lb(q) \in \boldsymbol{\Lambda}[q]$.
It will immediately be shown that 
$\ov{P}_\lb(q)$ is a well-defined element of 
$\boldsymbol{\Lambda}[q]$, i.e., $t_\lb^{(\Pi_s)}(q)$ divides $j(e^\lb{\cdot}\bD^{(s)}_q)$ in
$\boldsymbol{\Lambda}[q]$.
We say that   $\ov{P}_\lb(q)$ is a {\it short Hall-Littlewood polynomial}. (For, if $\Delta^+_s=
\Delta^+$ or if $\Delta^+_s$ and $\Pi_s$ are replaced with $\Delta^+$ and $\Pi$ 
in the above definition,
then one obtains the  usual Hall-Littlewood polynomials $P_\lb(q)$ for $\Delta$.)

\begin{prop}   \label{prop:P-delimost}
\[
 \ov{P}_\lb(q)=J\bigl(e^{\lb+\rho} \prod_{\ap\in\Delta^+_s,\,(\ap,\lb)>0} (1-qe^\ap)\bigr) J(\rho)^{-1} .
\]
\end{prop}\begin{proof}
1)  First consider the case in which $\lb=0$. Here
\[
 J(e^\rho){\cdot} j(e^0 {\cdot}\bD^{(s)}_q)=J\bigl(\sum_{A\subset \Delta^+_s} (-q)^{\# A}e^{\rho-|A|}\bigr).
\]
It is known that $\rho-|A|$ is regular if and only if $A=\mathsf{N}(w)$ for some $w\in W$ \cite{macd}.
Since $A\subset \Delta^+_s$, Lemma~\ref{lm:semi-dir}(ii) shows that actually $w\in W(\Pi_s)$.
Hence
\[
   J\bigl(\sum_{A\subset \Delta^+_s} (-q)^{\# A}e^{\rho-|A|}\bigr)=\sum_{w\in W(\Pi_s)}
   (-q)^{\ell(w)} J(e^{w^{-1}\rho})=\sum_{w\in W(\Pi_s)} q^{\ell(w)}{\cdot} J(e^\rho)=
   t_0^{(\Pi_s)}(q)J(e^\rho) .
\]
This proves that $\ov{P}_0(q)=1$.

2) \ For an arbitrary $\lb\in\Xp$, we  notice that 
$\sum_{w\in W_\lb} \esi(w) w(e^{\lb+\rho} \prod_{\ap\in\Delta^+_s}(1-qe^{-\ap}))$ is divisible by
$t_\lb^{(\Pi_s)}(q)$, by the first part of proof.
\\
(One has to consider the splitting $\prod_{\ap\in\Delta^+_s}(1-qe^{-\ap})=\prod_{\ap:\, (\ap,\lb)=0}(\dots) \prod_{\ap:\, (\ap,\lb)>0}(\dots)$, and use the fact
that $w(\prod_{\ap:\, (\ap,\lb)>0}(1-qe^{-\ap}))=\prod_{\ap:\, (\ap,\lb)>0}(1-qe^{-\ap})$
for any $w\in W_\lb$.)

This is already sufficient to conclude
that $\ov{P}_\lb(q)$ belongs to $\boldsymbol{\Lambda}[q]$. Further easy calculations that require a splitting $W\simeq W^\lb \times W_\lb$ 
are left to the reader.
\end{proof}

\begin{rmk}
Our proof is inspired by the remark in \cite[p.\,70, last paragraph]{rkg87}, where 
R.~Gupta refers to Macdonald's argument for the Hall-Littlewood symmetric functions.
\end{rmk}

\begin{rmk}
The Hall-Littlewood polynomials $P_\lb(q)$ 
interpolate between the irreducible characters $\chi_\lb$ (if $q=0$) and orbital sums
$\frac{1}{ \#(W_\lb)}\sum_{w\in W}e^{w\lb}$  (if $q=1$). For the short 
Hall-Littlewood polynomials $\ov{P}_\lb(q)$, we still have $\ov{P}_\lb(0)=\chi_\lb$.
At $q=1$, we obtain a linear combination of irreducible characters for $H$.
Namely, if $\chi^{(H)}_\mu$ denote the  character of $V^{(H)}_\mu$,
$\mu\in\mathfrak X_{+,H}$, then 
\[
 \ov{P}_\lb(1)= \frac{1}{ \#(W(\Pi_s)_\lb)} \sum_{w\in W(\Pi_s)} \chi^{(H)}_{w\lb} .
\]
An easy proof uses the semi-direct product structure of $W$ (Lemma~\ref{lm:semi-dir})
and Weyl's character formula for $H$.
(Note that if $\lb\in\Xp$, then $w\lb\in\mathfrak X_{+,H}$ for any $w\in W(\Pi_s)$.)
\end{rmk}

\begin{thm}    \label{thm:orthog-rel}
In $\boldsymbol{\Lambda}[[q]]$, the following relations hold:
\begin{enumerate}
\item   $\langle \ov{E}_\mu(q),\ov{P}_\lb(q)\rangle =\delta_{\lb,\mu}$; 
\item  $\ov{E}_\mu(q)=\displaystyle\frac{t_\mu^{(\Pi_s)}(q)}{\prod_{\ap\in\Delta_s}(1-qe^\ap)}
\ov{P}_\mu(q)$ and $\ov{E}_0(q)=\displaystyle\frac{t_0^{(\Pi_s)}(q)}{\prod_{\ap\in\Delta_s}(1-qe^\ap)}$.
\end{enumerate}
\end{thm}
\begin{proof}
(1) \ We mimic Gupta's proof of \cite[Theorem\,2.5]{rkg87}. The plan is as follows:
 \begin{itemize}
\item[\sf(i)] \  If $\chi_\pi$ occurs in $\ov{E}_\mu(q)=j(e^\mu{\cdot}\tilde\bD^{(s)}_q)$, then 
$\pi\curge \mu$;  and the coefficient of
$\chi_\mu$ equals $1$;

\item[\sf(ii)] \  If $\chi_\pi$ occurs in $j(e^\lb{\cdot}\bD^{(s)}_q)$, then $\pi\curle \lb$; and the coefficient of
$\chi_\lb$ equals $t_\lb^{(\Pi_s)}(q)$;

\item[\sf(iii)] \  Put $c_{\lb,\mu}=\langle  j(e^\lb{\cdot}\bD^{(s)}_q), j(e^\mu{\cdot}\tilde\bD^{(s)}_q)\rangle$.
Then $c_{\lb,\mu}=c_{\mu,\lb}$ and hence\\
$t_\mu^{(\Pi_s)}(q){\cdot}\langle \ov{E}_\lb(q),\ov{P}_\mu(q)\rangle= 
t_\lb^{(\Pi_s)}(q){\cdot}\langle \ov{E}_\mu(q),\ov{P}_\lb(q)\rangle.$
\end{itemize}
It will then follow that $c_{\lb,\mu}=\delta_{\lb,\mu}{\cdot}t_\lb^{(\Pi_s)}(q)$ proving the assertion.

For (i): By Weyl's character formula, the coefficient of $\chi_\pi$ in $j(e^\mu{\cdot}\tilde\bD^{(s)}_q)$
equals the coefficient of $e^{\pi+\rho}$ in (the expansion of) 
\[
J(e^\rho)\ov{E}_\mu(q) =\sum_{w\in W} \esi(w)w \left(\frac{e^{\mu+\rho}}{\prod_{\ap\in\Delta^+_s}
(1-qe^\ap)}\right) .
\]
This coefficient equals $\sum_{w,B} \esi(w)q^{\# B}$, where the summation is over
$w\in W$ and multi-sets $B$ of $\Delta^+_s$ such that $\pi+\rho=w(\mu+\rho+|B|)$.
Then $\pi+\rho \curge w^{-1}(\pi+\rho)=\mu+\rho+|B|\curge \mu+\rho$.
Hence $\pi\curge \mu$. If $\pi=\mu$, then the only possibility is $w=1$ and $B=\varnothing$.

For (ii):  Now, we are interested in the coefficient of $e^{\pi+\rho}$ in 
\[
  \sum_{w\in W} \esi(w) w\bigl(e^{\lb+\rho} \prod_{\ap\in\Delta^+_s}(1-qe^{-\ap})\bigr)
\]
It is equal to $\sum_{w, A} \esi(w)(-q)^{\# A}$, where the summation is over
$w\in W$ and subsets $A\subset\Delta^+_s$ such that $\pi+\rho=w(\lb+\rho-|A|)$.
Since $w\lb\curle \lb$ and $w(\rho-|A|)\curle \rho$, we obtain $\pi+\rho\curle \lb+\rho$.
Moreover, in case of equality we have $w\lb=\lb$ and $\rho-w^{-1}\rho=|A|$. This means
that $w\in W_\lb$ and $\mathsf{N}(w)=A\subset \Delta^+_s$. By Lemma~\ref{lm:semi-dir}(ii), we conclude that 
$w\in W(\Pi_s)$. Thus, $\#A=\ell(w)$ and 
the coefficient of $e^{\lb+\rho}$ equals
$\sum_{w\in W(\Pi_s)_\lb}   q^{\ell(w)}=t_\lb^{(\Pi_s)}(q)$.

For (iii): Set $\displaystyle\ov{\xi}=\frac{1}{\prod_{\ap\in\Delta_s}(1-qe^\ap)}$. 
It is a $W$-invariant element of $\boldsymbol{\Lambda}[[q]]$ and 
 $\bD^{(s)}_q \ov{\xi}=\tilde\bD^{(s)}_q$. Hence 
 $j(e^\mu{\cdot}\bD^{(s)}_q) \ov{\xi}=j(e^\mu{\cdot}\tilde\bD^{(s)}_q)$.
 But $\ov{\xi}$ is also a self-dual character. Thus, we have
 \[
   c_{\lb,\mu}=\langle j(e^\lb{\cdot}\bD^{(s)}_q), j(e^\mu{\cdot}\bD^{(s)}_q)\ov{\xi} \rangle=
   \langle j(e^\lb{\cdot}\bD^{(s)}_q)\ov{\xi}, j(e^\mu{\cdot}\bD^{(s)}_q) \rangle=c_{\mu,\lb} .
 \]
 (2) \ The equality 
 $\ov{E}_\mu(q)=t_\mu^{(\Pi_s)}(q)\ov{\xi}{\cdot}
\ov{P}_\mu(q)$ is essentially proved in (iii).  
Taking $\mu=0$ yields the rest.
\end{proof}

\begin{prop}    \label{prop:E-mu}
$\ov{E}_\mu(q)=\sum_{\lb\in\Xp} \ov{\me}_\lb^\mu(q) \chi_\lb$.
\end{prop}\begin{proof}
By definition,  
$\displaystyle J(e^\rho)\ov{E}_\mu(q)=J\left(\frac{e^{\mu+\rho}}{\prod_{\ap\in\Delta^+_s}
(1-qe^\ap)}\right)
=\sum_{\nu}\ov{\eus P}_q(\nu)J(e^{\mu+\nu+\rho})$. \\
The weight $\mu+\nu+\rho$ contributes to the last sum if and only if
$\mu+\nu+\rho=w(\lb+\rho)$ for some $\lb\in\Xp$ and $w\in W$. Hence
\begin{multline*}
  \sum_{\nu}\ov{\eus P}_q(\nu)J(e^{\mu+\nu+\rho})=
  \sum_{\lb\in\Xp}\sum_{w\in W} \ov{\eus P}_q(w(\lb+\rho)-(\mu+\rho))J(e^{w(\lb+\rho)})\\
  =\sum_{\lb\in\Xp}\sum_{w\in W} \esi(w)\ov{\eus P}_q(w(\lb+\rho)-(\mu+\rho))J(e^{\lb+\rho})=\sum_{\lb\in\Xp} \ov{\me}_\lb^\mu(q)J(e^{\lb+\rho}).
\end{multline*}
\end{proof}
\noindent
Part 1(ii) in the proof of Theorem~\ref{thm:orthog-rel} shows that $\{\ov{P}_\lb(q)\}_{\lb\in\Xp}$
is a $\BZ$-basis in $\boldsymbol{\Lambda}[q]$. Furthermore, 
Theorem~\ref{thm:orthog-rel}(1) and Proposition~\ref{prop:E-mu} readily 
imply that
\begin{equation}    \label{eq:xi-and-P}
         \chi_\pi=\sum_{\lb\in\Xp} \ov{\me}_\pi^\lb(q)\ov{P}_\lb(q) .
\end{equation}

Note that this sum is finite, since $\ov{\me}_\pi^\lb(q)=0$ unless $\lb \curle \pi$.
Let us transform the expression for $\ov{P}_\lb(q)$ given by definition:
\begin{multline*}
J(e^\rho){\cdot} t_\lb^{(\Pi_s)}(q) {\cdot}\ov{P}_\lb(q)=
J\bigl(e^{\lb+\rho}\prod_{\ap\in\Delta^+_s}(1-qe^{-\ap})\bigr)
\\
=J\left(e^\lb \frac{\prod_{\ap\in\Delta^+_s}(1-qe^{-\ap})}{\prod_{\ap>0} (1-e^{-\ap})}\cdot
e^\rho \prod_{\ap>0} (1-e^{-\ap})
\right)=\sum_{w\in W} w\left( e^\lb \frac{\prod_{\ap\in\Delta^+_s}(1-qe^{-\ap})}{\prod_{\ap>0} (1-e^{-\ap})} \right){\cdot} J(e^\rho) .
\end{multline*}
Hence 
$\ov{P}_\lb(q)= \displaystyle\frac{1}{t_\lb^{(\Pi_s)}(q)} \sum_{w\in W} w\left( e^\lb \frac{\prod_{\ap\in\Delta^+_s}(1-qe^{-\ap})}{\prod_{\ap>0} (1-e^{-\ap})} \right)$, and substituting this in Equation~\eqref{eq:xi-and-P} we obtain a generalisation of an  identity of Kato (cf. \cite[Theorem\,3.9]{rkg87}):
\begin{equation}    \label{eq:ttk}
   \chi_\pi=\sum_{\lb\in\Xp} \ov{\me}_\pi^\lb(q)\frac{1}{t_\lb^{(\Pi_s)}(q)} \sum_{w\in W} w
   \left( e^\lb \frac{\prod_{\ap\in\Delta^+_s}(1-qe^{-\ap})}{\prod_{\ap>0} (1-e^{-\ap})} \right) .
\end{equation}
Taking $q=1$, we obtain
\[
  \chi_\pi=\sum_{\lb\in\Xp}  \ov{m}_\pi^\lb {\cdot}\frac{1}{\# W(\Pi_s)_\lb} {\cdot} \sum_{w\in W} w
  \left( \frac{e^\lb}{ \prod_{\ap\in \Delta^+_l} (1-e^{-\ap}) }\right) .
\]
Taking into account that $W=W(\Pi_s)\rtimes W_l$ and 
$\ov{m}_\pi^\lb=\ov{m}_\pi^{w\lb}$ for any $w\in W(\Pi_s)$, this specialisation is equivalent
to the formula  $\chi_\pi=\sum_{\lb\in\mathfrak X_{+,H}}\ov{m}_\pi^\lb \chi_\lb^{(H)}$.

We introduce another bilinear form in $\boldsymbol{\Lambda}[q]$ such that
$\{\ov{P}_\lb(q)\}$ to be an orthogonal basis. 
To this end, the null-cone in $\vbt$ plays the same role as the nilpotent cone 
$\fN\subset\g$ for the Hall-Littlewood polynomials $P_\lb(q)$, cf. \cite[\S\,2]{rkg87-A}. 

For a graded $G$-module $\cM=\oplus_i \cM_i$ with $\dim\cM_i< \infty$, the graded character of $\cM$, $\mathsf{ch}_q(\cM)$, 
is the formal sum $\sum_i \mathsf{ch}(\cM_i) q^i \in \boldsymbol{\Lambda}[[q]]$.

\begin{prop}   \label{null-cone-little}
The graded character  
of the graded $G$-algebra $\BC[ \fN(\vbt)] $ equals 
 \[
   \mathsf{ch}_q(\BC[ \fN(\vbt)])=\frac{t_0^{(\Pi_s)}(q)}{\prod_{\ap\in\Delta_s}(1-qe^\ap)}=t_0^{(\Pi_s)}(q) {\cdot}\ov{\xi}= \ov{E}_0(q) .
\]
 \end{prop}
 \begin{proof}  The weight structure of $\vbt$ (Lemma~\ref{lm:vesa-little}) shows that the graded character of $\BC[\vbt]$ equals 
 $\displaystyle \mathsf{ch}_q(\BC[ \vbt])=\frac{1}{(1-q)^{\#\Pi_s} \prod_{\ap\in\Delta_s}(1-qe^\ap)}$. We know that $\fN(\vbt)$ is a  complete intersection
 of codimension $m:=\# \Pi_s$ and the ideal of $\fN(\vbt)$ is generated by 
 algebraically independent generators of $\BC[\vbt]^G$. Furthermore,
 if $d_1,\dots, d_m$  are the degrees of these generators, then 
 $d_1-1,\dots, d_m-1$ are the exponents of $W(\Pi_s)$ (Prop.~\ref{nullcone-vbt}). Thus,
 \[
 \mathsf{ch}_q(\BC[ \fN(\vbt)])=\frac{\prod_{i=1}^m(1-q^{d_i})}{(1-q)^{m} \prod_{\ap\in\Delta_s}(1-qe^\ap)}=
 \frac{\prod_{i=1}^m(1+q+\dots +q^{d_i-1})}{\prod_{\ap\in\Delta_s}(1-qe^\ap)},
 \]
 and it is well known that  $t^{(\Pi_s)}_0(q)=\prod_{i=1}^m(1+q+\dots +q^{d_i-1})$.
 \end{proof}  
 
Combining Propositions~\ref{prop:E-mu} and \ref{null-cone-little} yields
\[
     \mathsf{ch}_q(\BC[ \fN(\vbt)])=\sum_{\lb\in\Xp} \ov{\me}_\lb^0(q) \chi_\lb ,
\]
which is \cite[Theorem\,4]{ww}.  In other words,  
$\sum_{i\ge 0}\dim \bigl(\text{Hom}_G(V_\lb, \BC[\fN(\vbt)]_i)\bigr)q^i= 
\ov{\me}_\lb^0(q)$ for every $\lb\in\Xp$.
 
 Define a new bilinear form in $\boldsymbol{\Lambda}[q]$ by letting 
 \[
 \langle\!\langle \chi_\lb,\chi_\mu\rangle\!\rangle= \langle \chi_\lb\chi_\mu^*,\, 
 t_0^{(\Pi_s)}(q) {\cdot}\ov{\xi} \rangle=\langle \chi_\lb,\, 
 t_0^{(\Pi_s)}(q) {\cdot}\ov{\xi}\chi_\mu \rangle .
 \]
 In view of Proposition~\ref{null-cone-little}, $ \langle\!\langle \chi_\lb,\chi_\mu\rangle\!\rangle$
 is a polynomial in $q$ that counts graded occurrences of the $G$-module $V_\lb\otimes
 V_\mu^*$ in $\BC[ \fN(\vbt)]$.
 
 \begin{thm}    \label{thm:new-scalar-P}
$\displaystyle
\langle\!\langle \ov{P}_\lb(q), \ov{P}_\mu(q) \rangle\!\rangle=\frac{t_0^{(\Pi_s)}(q)}{t_\mu^{(\Pi_s)}(q)}\delta_{\lb,\mu}$.
\end{thm}\begin{proof}
By definition and Theorem~\ref{thm:orthog-rel}, we have
\[
  \langle\!\langle \ov{P}_\lb(q), \ov{P}_\mu(q) \rangle\!\rangle=
  \langle \ov{P}_\lb(q), t_0^{(\Pi_s)}(q) {\cdot}\ov{\xi}{\cdot}\ov{P}_\mu(q) \rangle =
  \langle \ov{P}_\lb(q),\frac{t_0^{(\Pi_s)}(q)}{t_\mu^{(\Pi_s)}(q)} \ov{E}_\mu(q) \rangle =
  \frac{t_0^{(\Pi_s)}(q)}{t_\mu^{(\Pi_s)}(q)}\delta_{\lb,\mu}.
\]
Here we also use the fact that $\bD^{(s)}_q \ov{\xi}=\tilde\bD^{(s)}_q$ and  hence
$t_\mu^{(\Pi_s)}(q){\cdot}\ov{\xi}{\cdot}\ov{P}_\mu(q)=\ov{E}_\mu(q)$.
\end{proof}

Finally, using Eq.~\eqref{eq:xi-and-P}, we obtain
\[
   \langle\!\langle \chi_\lb,\chi_\mu\rangle\!\rangle=\sum_{\pi\in\Xp}
   \ov{\me}_\lb^\pi(q)\, \ov{\me}_\mu^\pi(q)\, \frac{t_0^{(\Pi_s)}(q)}{t_\pi^{(\Pi_s)}(q)} .
\]

\section{Miscellaneous remarks } 
\label{sect:miscell}

\subsection{}
It is noticed in \cite[5.1]{rkg87} that Lusztig's $q$-analogues $\me_{\lb}^\mu(q)$ satisfy the identity
\begin{equation}    \label{eq:series-Lqa}
   \sum_{\mu\in\mathfrak X} \me_{\lb}^\mu(q) e^\mu=\frac{J(e^{\lb+\rho})}{e^\rho\prod_{\ap>0}(1-qe^{-\ap})}=
   \chi_\lb\cdot \prod_{\ap>0} \frac{(1-e^{-\ap})}{(1-qe^{-\ap})} .
\end{equation}
This can be regarded as quantisation of the equality $\chi_\lb=\sum_\mu m_\lb^\mu e^\mu$,
which describes $V_\lb$ as $T$-module.
In the context of short $q$-analogues, we wish to have a quantisation of the
equality $\chi_\lb=\sum_{\mu\in \mathfrak X_{+,H}} \ov{m}_\lb^\mu \chi^{(H)}_\mu$, which describes 
$V_\lb$ as $H$-module \cite[\S\,3]{heck}.  
The desired quantisation is
\begin{prop}
   $\displaystyle
   \sum_{\mu\in\mathfrak X_{+,H}}\ov{\me}_\lb^\mu(q) \chi^{(H)}_\mu=\chi_\lb{\cdot}
   \frac{1}{\# W_l}
   \sum_{w\in W_l} w\left( \prod_{\ap\in\Delta^+_s}\frac{1-e^{-\ap}}{1-qe^{-\ap} } \right)$ .
\end{prop}
\begin{proof}  Using Weyl's formula, the function 
$(\mu\in \mathfrak X_{+,H}) \mapsto \chi^{(H)}_\mu$ can be extended to the whole of $\mathfrak X$ such that it will satisfy the identity
$ \chi^{(H)}_{w\odot\mu}=\esi(w) \chi^{(H)}_\mu$, $w\in W_l$. Recall that `$\odot$' stands for the
shifted action of $W_l$. Since the same identity holds for $\ov{\me}_\lb^\mu(q)$,
see Eq.~\eqref{eq:shifted_act}, the left hand side can be replaced with
$\displaystyle
\frac{1}{\# W_l}   \sum_{\mu\in\mathfrak X}\ov{\me}_\lb^\mu(q) \chi^{(H)}_\mu$.
The rest can by achieved via routine transformations of this sum, 
using the definition of $\ov{\me}_\lb^\mu(q)$ 
and  Weyl's character formulae for $H$ and $G$.
\end{proof}
Yet another quantisation, which is easier to prove, is 
\begin{equation}    \label{eq:series-sqa}
    \sum_{\mu\in\mathfrak X} \ov{\me}_{\lb}^\mu(q) e^\mu=
    \frac{J(e^{\lb+\rho})}{e^\rho\prod_{\ap\in\Delta^+_s}(1-qe^{-\ap})}=
    \chi_\lb\cdot \frac{\prod_{\ap>0} (1-e^{-\ap})}{\prod_{\ap\in\Delta^+_s}(1-qe^{-\ap})} .
\end{equation}
Comparing Equations~\eqref{eq:series-Lqa} and \eqref{eq:series-sqa}, we obtain a relation
between  Lusztig's and short $q$-analogues:
\[
   \prod_{\ap\in\Delta^+_l}(1-qe^{-\ap}) \sum_{\mu\in\mathfrak X} \me_{\lb}^\mu(q) e^\mu=
   \sum_{\nu\in\mathfrak X} \ov{\me}_{\lb}^\nu(q) e^\nu .
\]
Whence  $\displaystyle \ov{\me}_{\lb}^\mu(q)=\sum_{A\subset \Delta^+_l} (-q)^{\# A} \me_{\lb}^{\mu+|A|}(q)$. Or, conversely,
$\displaystyle {\me}_{\lb}^\mu(q)=\sum_{B} q^{\# B} \ov{\me}_{\lb}^{\mu+|B|}(q)$, where $B$ ranges over the 
finite multisets in $\Delta^+_l$.  In particular, taking $q=1$ and $\mu=0$, we  obtain
\[
     \dim V_\lb^H = \ov{m}_\lb^0=\sum_{A\subset \Delta^+_l} (-1)^{\# A} m_{\lb}^{|A|} .
\]
\begin{ex-bn}
If $G=Sp_{2n}$, then $H=(SL_2)^n$ and $\Delta^+_l=\{2\esi_1,\dots,2\esi_n\}$.
Here $\esi_{i_1}+\ldots +\esi_{i_k}$ is $W$-conjugate to  $\vp_k=\esi_1+\ldots +\esi_k$ and the previous relation becomes
\[
    \dim V_\lb^H = \sum_{k=0}^n (-1)^k \genfrac{(}{)}{0pt}{}{n}{k} m_\lb^{2\vp_k} . 
\]
\end{ex-bn}

\subsection{}
It is well known that, for $\lb$ strictly dominant, the Hall-Littlewood polynomials $P_\lb(q)$ have a nice specialisation at $q=-1$:
If $\lb\gtrdot\rho$, then $P_\lb(-1)=\chi_{\lb-\rho} \chi_\rho$.  (See \cite[7.4]{visw} for a generalisation to symmetrisable Kac-Moody algebras.)
For $\Delta$ of type $\GR{A}{n}$, $P_\lb(-1)$ is a classical Schur's $Q$-function
\cite[III.8]{macdo95}.
 A similar phenomenon occurs for short Hall-Littlewood polynomials.
 
\begin{prop}    \label{short_HL_at_q=-1}
Suppose $\lb\gtrdot \rho_s$ and $G$ is of type $\GR{B}{n},\GR{C}{n}$, or $\GR{F}{4}$.
Then $\ov{P}_\lb(-1)= \chi_{\lb-\rho_s}\chi_{\rho_s}$.
\end{prop}
\begin{proof}
If $\lb\gtrdot \rho_s$, then $t^{(\Pi_s)}_\lb(q)=1$ and 
\begin{multline*}
\ov{P}_\lb(-1)=
J\bigl(e^{\lb+\rho}\prod_{\lb\in\Delta^+_s}(1+e^{-\ap})\bigr) J(e^\rho)^{-1}=
\\
\sum_{w\in W} \esi(w) w(e^{\lb-\rho_s+\rho})(e^{\rho_s}\prod_{\lb\in\Delta^+_s}(1+e^{-\ap}))
{\cdot}J(e^\rho)^{-1}= \chi_{\lb-\rho_s}{\cdot}\prod_{\lb\in\Delta^+_s}(e^{\ap/2}+e^{-\ap/2}) .
\end{multline*}
For $G$ is of type $\GR{B}{n},\GR{C}{n}$, or $\GR{F}{4}$, it is known that 
$\chi_{\rho_s}=\prod_{\lb\in\Delta^+_s}(e^{\ap/2}+e^{-\ap/2})$ \cite[Theorem\,2.9]{spin}.
\end{proof}

\begin{rmk}
The  proof of equality $\chi_{\rho_s}=\prod_{\lb\in\Delta^+_s}(e^{\ap/2}+e^{-\ap/2})$ in
\cite{spin} is only based on the assumption that $\|\text{long}\|^2/\|\text{short}\|^2=2$, 
i.e., it does not refer to classification. 
For $\GR{G}{2}$, the true equality is $\prod_{\lb\in\Delta^+_s}(e^{\ap/2}+e^{-\ap/2})=
\chi_{\rho_s}+1$.
\end{rmk}

\subsection{}  \label{subs:jump}
Ranee Brylinski proved that Lusztig's $q$-analogues $m_\lb^\mu(q)$ can be
computed via a principal filtration on $V_\lb^\mu$ whenever 
$H^i(G\times_B \ut, \cl_{G\times_B \ut}(\BC_\mu)^\star)=0$ for all $i\ge 1$. 
Namely, $\me_\lb^\mu(q)$ coincides with the ``jump polynomial" of the
principal filtration, see \cite{rkb89} for details. Another approach to her results can be found
in \cite{JLZ}.

I hope that a similar description exists for short $q$-analogues. First, we need a subspace of
$V_\lb$ whose dimension equals $\ov{m}_\lb^\mu=\mathsf{mult}( V^{(H)}_\mu, V_\lb)$.
Let  $V_\lb^{U(H)}$ be the subspace of $H$-highest vectors in $V_\lb$ with respect to
$\Delta^+_l$.  Then $V_\lb^{U(H),\mu}=V_\lb^{U(H)}\cap V_\lb^\mu$ has the required dimension.  For $\ap\in\Delta^+$, let $e_\ap$ be a nonzero root vector of $\g$.
Brylinski's principal filtration is determined by the principal nilpotent element
$e=\sum_{\ap\in\Pi}e_\ap$. In the context of short $q$-analogues, we consider 
$e_s=\sum_{\ap\in\Pi_s}e_\ap$ and the corresponding filtration of $V_\lb^{U(H),\mu}$.
That is, we set 
\[
    J^p_{e_s}(V_\lb^{U(H),\mu})=\{ v\in V_\lb^{U(H),\mu}\mid e_s^{p+1}{\cdot}v=0\} .
\]
The jump polynomial is defined to be 
\[
    \ov{r}_\lb^\mu(q)= \sum_{p\ge 0} 
    \dim\bigl( J^p_{e_s}(V_\lb^{U(H),\mu})/J^{p-1}_{e_s}(V_\lb^{U(H),\mu})\bigr)q^p .
\]
\begin{conj}
If $\mu\in\mathfrak X_{+,H}$ satisfies vanishing conditions of Theorem~\ref{thm:main}, then
$\ov{r}_\lb^\mu(q)=\ov{\me}_\lb^\mu(q)$.
\end{conj}

\subsection{}  Although the collapsing $f: \bZ=G\times_B \vbt^+\to \fN(\vbt)$  is not generically finite, it can be used for deriving useful properties of the null-cone.
Let $\varrho: \co_{\fN(\vbt)} \to Rf_\ast \co_{\bZ}$ be the corresponding natural morphism.
Since $f$ is projective, $H^0(\bZ,\co_\bZ)$ is a finite $\BC[\fN(\vbt)]$-module; 
and there is the trace map $H^0(\bZ,\co_\bZ)\to \BC[\fN(\vbt)]$ 
because $\fN(\vbt)$ is normal.  The trace map determines a morphism 
(in the derived category of $\co_{\fN(\vbt)}$-modules)
$\varrho':  Rf_\ast \co_{\bZ}\to \co_{\fN(\vbt)}$.
By Theorem~\ref{thm:main}, $H^i(\bZ,\co_\bZ)=0$ for $i\ge 1$, i.e., 
$R^i f_\ast\co_\bZ=0$ for $i\ge 1$.
Hence 
$\varrho'\circ\varrho$ is a quasi-isomorphism of $\co_{\fN(\vbt)}$ with itself.
Therefore, by 
\cite[Theorem\,1]{kov00},  $\fN(\vbt)$ has only rational singularities.

Clearly, this argument works in a more general context and yields the following:

\begin{prop} Let $N$ be a $P$-stable subspace in a $G$-module $V$.
If  $H^i(G\times_P N,\co_{G\times_P N})=0$
for all $i\ge 1$, then the normalisation of $G{\cdot}N$ has only rational singularities. 
\end{prop}

\end{document}